\documentclass[11pt,twoside]{article}
\usepackage{amsmath,bm}
\usepackage{amssymb}
\usepackage{fancyhdr}
\usepackage{latexsym}
\usepackage{bbding}
\usepackage{mathrsfs}
\usepackage{exscale}
\usepackage{relsize}
\usepackage{wasysym}
\usepackage{cite}
\usepackage{multicol,graphics}
\usepackage{color}
\allowdisplaybreaks[4]
\def\to{\rightarrow}

\def\fr{\frac}

\def\i{\infty}

\def\R{\mathbb{R}}

\def\p{\partial}

\def\lb{\label}

\def\bar{\overline}

\providecommand{\Nop}{\mathscr N}
\providecommand{\supp}{\operatorname{supp}}
\providecommand{\diver}{\operatorname{div}}
\providecommand{\ess}{\operatorname{ess}}

\newcommand{\SE}{\setcounter{equation}{0} \section}
\renewcommand\tilde{\widetilde}
\newcommand{\be}{\begin{equation}}
\newcommand{\ee}{\end{equation}}
\newcommand{\bc}{\begin{cases}}
\newcommand{\ec}{\end{cases}}
\newcommand{\bes}{\begin{equation*}}
\newcommand{\ees}{\end{equation*}}
\newcommand{\bls}{\begin{align*}}
\newcommand{\els}{\end{align*}}
\newcommand{\baa}{\begin{array}}
\newcommand{\eaa}{\end{array}}
\newcommand{\ba}{\begin{eqnarray}}
\newcommand{\ea}{\end{eqnarray}}
\newcommand{\bas}{\begin{eqnarray*}}
\newcommand{\eas}{\end{eqnarray*}}
\newcommand{\bd}{\begin{description}}
\newcommand{\ed}{\end{description}}

\newtheorem{theo}{\bf Theorem}[section]
\newtheorem{lem}[theo]{\bf Lemma}
\newtheorem{prop}[theo]{\bf Proposition}

\newtheorem{defi}[theo]{\bf Definition}
\newtheorem{rem}[theo]{\bf Remark}

\newenvironment{pr}[1][Proof]{\noindent\textbf{#1.} }{\hfill $\Box$}
\allowdisplaybreaks
\begin{document}
\date{}
\title{\bf{Propagation phenomena in spatial periodic bistable
reaction-diffusion equations around an obstacle}}
\author{Wei-Jie Sheng\thanks{Corresponding author
 (E-mail address: shengwj09@hit.edu.cn).}
\\
\footnotesize{ School of Mathematics, Harbin Institute of Technology}, \\
\footnotesize{Harbin, Heilongjiang, 150001, People's Republic of China}
}
\maketitle

\begin{abstract}
\noindent

This paper is concerned with propagation phenomena for a spatially periodic
bistable reaction--diffusion equation around an obstacle. It is known that an
entire solution emanating from a pulsating front converges to the same front
as $t\to+\i$ under complete propagation and an additional
technical assumption. We remove the technical assumption and prove that 
the entire solution coming from a pulsating front converges to the same pulsating
front as \(t\to+\infty\) under the complete propagation.
\end{abstract}

\section{Introduction and main results}

This paper is concerned with the propagation phenomena for the following
reaction--diffusion problem:
\be\lb{eq1.1}
\bc
u_t=\Delta u+f(x,u)\ \ \text{for $t\in\R$ and $x\in\Omega$},\\
\nu\cdot\nabla u=0\ \ \text{for $t\in\R$ and $x\in \p\Omega=\p K$},
\ec
\ee
where $\Omega=\R^N\backslash K\subset\R^N$ is  an exterior domain,
$K$ is a compact set which is the closure of an open set
with smooth boundary,
and the Neumann no flux boundary condition
is imposed, $\nu=\nu(x)$ is the outward unit normal on the boundary $\p \Omega$.  
 Throughout, we work in dimension $N\geq2$.
We assume that the nonlinearity $f$ satisfies the following conditions.
\begin{description}
 \item[(f1)] (smooth) \( f(x, u) \in C(\mathbb{R}^N \times [0, 1]; \mathbb{R}) \) is of class \( 
    C^\alpha \) with respect to spatial variable \( x \) uniformly for \( u \in [0, 1] \) with \( \alpha \in 
    (0, 1) \), and of class \( C^1 \) with respect to variable \( u \in [0, 1] \) uniformly in \( x \in 
    \mathbb{R}^N \), while \( f_u(x, u) = \partial_u f(x, u) \) is Lipschitz continuous with respect to 
    \( u \in [0, 1] \) uniformly for \( x \in \mathbb{R}^N \).
\item[(f2)] (\( L \)-periodic) There is a vector \( L := (L_1, \ldots, L_N) \) with \( L_i > 0 \) (\( i 
    = 1, \ldots, N \)) such that for any vector \( k := (k_1, \ldots, k_N) \in \mathbb{Z}^N \),
    \[
    f(x + kL, u) = f(x, u),
    \]
    where \( kL := (k_1 L_1, \ldots, k_N L_N) \).
\item[(f3)] (bistable) \( f(x, 0) = f(x, 1) = 0 \) and there are \( \eta_0 \in (0, 1/4) \) and \( \beta_0 > 0 
\)   such that
 \be\lb{eq1.2}
-f_u(x, u) \geq \beta_0/2 \quad \text{for all} \ (x, u) \in \mathbb{R}^N \times ([0, 2\eta_0] \cup [1 
 -  2\eta_0, 1]).
 \ee
\end{description}
In particular,
\[
\max\{f_u(x,0),f_u(x,1)\}\leq-\beta_0/2,
\ \  x\in\mathbb R^N.
\]
For mathematical convenience, the function \(f\) is extended in \(\Omega\times(\R 
\backslash[0,1])\) as 
follows:
\be\lb{eq1.6}
\begin{array}{ll}
f(x, u)=\partial_{u} f(x, 0) u & \text { for }(x, u) \in \Omega \times(-\infty, 0), \\
f(x, u)=\partial_{u} f(x, 1)(u-1) & \text { for }(x, u) \in \Omega \times(1,+\infty) .
\end{array}
\ee
Thus \(f\) is continuous in \(\Omega \times \mathbb{R}\), $L$-periodic in \(x\), \(\min _{x \in 
\Omega} f(x, u)>0\) for all \(u<0\), and \(\max _{x \in \Omega} f(x, u)<0\) for all \(u>1\), 
A particular example of such a function \(f\) satisfying (f1)-(f3) is the cubic nonlinearity
\begin{equation*}
f(x, u)=u(1-u)\left(u-\theta_{x}\right)
\end{equation*}
where \(x \mapsto \theta_{x}\) is a 1 -periodic \(C^{0, \alpha}(\Omega)\) function ranging in 
\((0,1)\).

Jia, Sheng and Wang \cite{jsw} established an entire
solution emanating from a pulsating front  $U_{e_1}(x,x_1-c_{e_1}t)$, whose definition is recalled
below. Under complete propagation and the additional condition
\be\lb{add}
\sup_{\substack{(x,t)\in\overline\Omega\times\mathbb R\\
x_1-c_{e_1}t\geq\mathcal R}}
\left|\frac{u(x,t)}{U_{e_1}(x,x_1-c_{e_1}t)}-1\right|
\longrightarrow0\quad\text{as }\mathcal R\to+\infty,
\ee
they proved uniform convergence to the same front as $t\to+\infty$.
Condition \eqref{add} controls the relative error far ahead of the front,
uniformly over all times. Shi and Li \cite{shi-li-stability} removed \eqref{add} when the reaction
term depends periodically on the propagation coordinate alone. Their
equation is posed in $\mathbb R^N$, with reaction term $f(x_N,u)$ independent
of the transverse variables. They proved stability under decaying
perturbations in the whole space and  obtained recovery
in exterior domains under complete propagation.
The extension to a general spatially periodic reaction term $f(x,u)$ was
left open by their method.

In this paper, we extend this recovery result to reaction terms that depend
periodically on all spatial coordinates. Under the standing assumptions
and complete propagation, we prove uniform convergence to the incoming
pulsating front without \eqref{add}. The proof combines far-field
convergence with a positive periodic adjoint weight and a weighted
dissipation estimate for localized phase perturbations. The comparison
principle then yields recovery with the phase prescribed by the incoming
front.

Let us first recall the notion of a pulsating traveling front for the
whole space equation
\be\lb{eq1.5}
u_t=\Delta u+f(x,u),\ \  x\in\mathbb R^N,\ t\in\mathbb R.
\ee
Write
$\mathbb T^N=\prod_{i=1}^N(\mathbb R/L_i\mathbb Z)$
for the torus associated with the spatial periods.

\begin{defi}
A moving pulsating traveling front of \eqref{eq1.5} in the direction
$e\in\mathbb S^{N-1}$ is a pair $(U_e,c_e)$, where $c_e\ne0$ and
$U_e:\mathbb T^N\times\mathbb R\to(0,1)$, such that
$u(x,t)=U_e(x,x\cdot e-c_et)$ is a classical entire solution of
\eqref{eq1.5}. The profile $U_e(x,s)$ is $L$-periodic in $x$ for each fixed
$s$ and satisfies
\be\lb{eq+1.5}
\lim_{s\to+\infty}U_e(x,s)=0\ \ \text{and}\ \
\lim_{s\to-\infty}U_e(x,s)=1
\quad\text{uniformly in }x\in\mathbb T^N.
\ee
\end{defi}
It is obvious that the profile satisfies
\[
c_e\partial_sU_e+\partial_{ss}U_e+\Delta_xU_e
+2e\cdot\nabla_x\partial_sU_e+f(x,U_e)=0,
\ \  x\in\mathbb T^N,\ s\in\mathbb R.
\]

The definition was first introduced by Shigesada et al. \cite{shi86},
Xin \cite{xin92a,xin92b,xin92c}, and developed by Berestycki and Hamel \cite{bh}.
Ducrot \cite{duc} proved that either there is a moving pulsating traveling front  $(U_e, 
c_e)$
in the direction $e$ with  ${\rm sign}(c_e)={\rm sign}(\mathcal{I})$, or there is an
almost monotonic  standing transition solution $U_e \equiv U_e(x)$ in the direction $e$ when 
$\mathcal I \neq  0$, where
\begin{equation*}
\mathcal{I}=\int_{\mathbb{T}^{N}} d x \int_{0}^{1} f(x, s) d s,
\end{equation*}
which further means that the unbalanced condition (namely $\mathcal I \neq  0$)
is not sufficient to ensure the existence of  a moving pulsating traveling wave solution with 
effective speed,
although it is true in the one-dimensional framework \cite{dhz1}. 
Under an additional assumption (see Assumption 1.3 in \cite{gr}),
Giletti and Rossi \cite{gr} proved that for any
$e \in \mathbb S^{N-1}$, there exists a monotonic in time pulsating traveling front connecting $1$ 
to $0$ in the direction $e$ with some speed $c_e\in\R$.
Since we are interested in the convergence of the solution coming from a pulsating front, hence we
impose the following assumptions.
\bd
\item[(A1)]
For every $e\in\mathbb S^{N-1}$, equation \eqref{eq1.5} admits a moving pulsating traveling front $U_e(x,x\cdot e-c_et)$ with speed $c_e\ne0$.
\item[(A2)] $\mathcal{I}>0$.
\ed
Under conditions (A1) and (A2), Guo \cite{guo} proved  that the speed $c_e$ and the profile
{$U_e(x,s)$} are continuous and differentiable with respect to 
the direction $e$, and also he gave an estimate of the speed 
\be\lb{signspeed}
 {0<\min_{e\in \mathbb S^{N-1}} c_e\leq
\max_{e\in \mathbb S^{N-1}} c_e<+\infty.}
\ee
Guo et al. \cite{gllw} further considered the existence, uniqueness and stability of curved fronts 
for bistable reaction-diffusion equations in spatially periodic media under the conditions (A1) and 
(A2). Here we restrict attention to entire solutions emanating from a single pulsating front.
For the sake of convenience, we just consider the direction \(e = e_1 = (1, 0,..., 0) \in \mathbb 
S^{N-1}\)  and always denote \(x_1 = x \cdot e_1\), $c_1=c_{e_1}$ and $U_1=U_{e_1}$ in the 
sequel. That is, we consider the pulsating traveling front $(U_1(x,x_1-c_1t),c_1)$ which satisfies
\begin{equation*}
c_1 \partial_{s} U_1+\p_{ss} U_1+\Delta_x U_1+2 \p_{x_1} (\p_{s}U_1)+f(x,U_1)=0, 
\quad \forall(s, x) \in \mathbb{R} \times \mathbb{T}^{N}.
\end{equation*}
We always remember that \(c_1 > 0\) by \eqref{signspeed}.

Call \(k_0(\lambda)\) and \(k_1(\lambda)\)  the principal eigenvalues of the operator
\[
L_{\lambda_0} \psi:=-\Delta \psi+2 \lambda e  \cdot\nabla \psi -\left[\lambda^{2}
+f_u(x,0)\right] \psi
\]
and 
\[
L_{\lambda_1} \psi:=-\Delta \psi+2 \lambda e  \cdot\nabla \psi -\left[\lambda^{2}
+f_u(x,1)\right] \psi,
\]
respectively, acting on the set 
\(E_{\lambda}=\left\{\psi \in C^{2}(\R^N), \psi  \ \text{is periodic in } \R^N  \right\}\).
It can be inferred from \cite[Lemma 2.1]{h} and \cite[Proposition 5.7]{bh} that 
$k_0(\lambda)$  and $k_1(\lambda)$  are
analytic and concave in $\R$, respectively,  and there hold
\[
\{\lambda\in\R:k_0(\lambda)-c\lambda=0\}=\{\lambda_0^-,\lambda_0^+\}
\]
and
\[
\{\lambda\in\R:k_1(\lambda)-c\lambda=0\}=\{\lambda_1^-,\lambda_1^+\}
\]
with $\lambda_0^-<0<\lambda_0^+$ and $\lambda_1^-<0<\lambda_1^+$. 
Let \(\psi_{\lambda_0^+}\) and 
\(\psi_{\lambda_1^+}\) denote the unique positive eigenfunctions
associated with the eigenvalues  $\lambda_0^+$ and $\lambda_1^+$, respectively,
such that 
\[
\left\|\psi_{\lambda_0^+}\right\|_{L^{\infty}(\mathbb T^N)}=1\ \ 
\text{and}\ \ \left\|\psi_{\lambda_1^+}\right\|_{L^{\infty}(\mathbb T^N)}=1.
\]
A similar treatment as that of \cite[Theorem 1.3]{h}
and an application of L' H\^{o}spital's rule
yields the following result.
\begin{prop}\lb{prop1}
There exist some positive constant $\rho_0$ and $\rho_1$ such that
\begin{equation*}
\lim_{s\to+\i}\fr{U_1(x,s)}{\rho_0 e^{-\lambda_0^+ s} \psi_{\lambda_0^+}(x)}=1 \ \ \text{and} \ \
\lim_{s\to-\i}\fr{1-U_1(x,s)}{\rho_1 e^{\lambda_1^+ s} \psi_{\lambda_1^+}(x)}=1\ \ 
\text{uniformly in $x\in\mathbb T^N$},
\end{equation*}
\begin{equation*}
\lim_{s\to+\i}\fr{\p_s U_1(x,s)}{\rho_0 \lambda_0^+ e^{-\lambda_0^+ s} \psi_{\lambda_0^+}(x)}=-1 \ \ \text{and} \ \
\lim_{s\to-\i}\fr{\p_s U_1(x,s)}{\rho_1 \lambda_1^+ e^{\lambda_1^+ s} 
\psi_{\lambda_1^+}(x)}=-1\ \ \text{uniformly in $x\in\mathbb T^N$},
\end{equation*}
\begin{equation*}
\lim_{s\to+\i}\fr{\p_{ss} U_1(x,s)}{\rho_0 (\lambda_0^+)^2 e^{-\lambda_0^+ s} \psi_{\lambda_0^+}(x)}=-1 \ \ \text{and} \ \
\lim_{s\to-\i}\fr{\p_{ss} U_1(x,s)}{\rho_1 (\lambda_1^+)^2 e^{\lambda_1^+ s} 
\psi_{\lambda_1^+}(x)}=1\ \ \text{uniformly in $x\in\mathbb T^N$}.
\end{equation*}
\end{prop}

Besides pulsating traveling fronts, an effective way to study propagation phenomena in general 
unbounded domains is transition fronts, since pulsating fronts may not exist in general unbounded
domains.
The notions of transition fronts and their global mean speed
are first introduced in the seminal works \cite{bh,bh1,bh2} by Berestycki and Hamel  in high-dimensions, see also Shen
\cite{s} for the one-dimensional setting. Before going further, let us first recall the notions of 
transition fronts
and their global mean speed.
For any two subsets $A$ and $B$ in $\overline\Omega$, we set
\begin{equation*}
d_\Omega(A,B)=\inf\{d_\Omega(x,y): (x, y)\in A\times B\}\mbox{ and }  
d_\Omega(x,A)={d_\Omega(\{x\},A)},
\end{equation*}
where $d_\Omega$ is the geodesic distance in $\overline\Omega$. Consider two families of 
nonempty unbounded open subsets  $\{\Omega^\pm_t\}_{t\in\R}\subset\mathbb R^N$ satisfying
\begin{equation}\label{subset1}
\forall t\in\mathbb R,\,\left\{\begin{array}{l}
\Omega_t^-\cap\Omega_t^+=\emptyset,\\
\partial\Omega_t^\pm\cap\Omega=:\Gamma_t,\\
\Omega_t^-\cup\Gamma_t\cup\Omega_t^+=\Omega,\\
\sup\{d_\Omega(x,\Gamma_t):x\in\Omega_t^-\}=\sup\{d_\Omega(x,\Gamma_t):x\in\Omega_t^+\}=\infty
\end{array}\right.
\end{equation}
and
\begin{equation}\label{subset2}
\left\{
\begin{array}{ll}
\inf\left\{\sup\{d_\Omega(y,\Gamma_t):y\in\Omega_t^+,d_\Omega(x,y)\leq 
r\}:t\in\R,x\in\Gamma_t\right\}\to\infty\ \mbox{ as }\ r\to\infty,\\
\inf\left\{\sup\{d_\Omega(y,\Gamma_t):y\in\Omega_t^-,d_\Omega(x,y)\leq 
r\}:t\in\R,x\in\Gamma_t\right\}\to\infty\ \mbox{ as }\ r\to\infty.
\end{array}\right.
\end{equation}
Notice that the condition \eqref{subset1} implies that $\Gamma_t$ divides the domain $\Omega$ 
into two disjoint unbounded open sets $\Omega_t^-$ and $\Omega_t^+$, which also ensures the 
interface $\Gamma_t$ is nonempty for each $t\in\mathbb R$. As far as \eqref{subset2} is 
concerned, it means that for any $M>0$, there exists a constant $r_M>0$ such that for each 
$t\in\mathbb R$ and $x\in\Gamma_t$, there are $y^\pm\in\Omega$ such that
\begin{align*}\label{subset2-2}
y^\pm\in\Omega^\pm_t,\,d_\Omega(x,y^\pm)\leq r_M\text{ and 
}d_\Omega(y^\pm,\Gamma_t)\geq M.
\end{align*}
In other words, for each $t\in\R$, any point on $\Gamma_t$ is not too far away from the centers of 
two large balls included in $\Omega_t^-$ and $\Omega_t^+$. Moreover, the sets $\Gamma_t$ are 
assumed to be included in a
finite number of graphs. Namely, there is an integer $n\geq1$ such that for each $t\in\mathbb R$, 
there exist $n$ open subsets $\omega_{i,t}\subset\R^{N-1}$ ($1\leq i\leq n$),
$n$ continuous maps $\psi_{i,t}:\omega_{i,t}\to\R$ and $n$ rotations $R_{i,t}$ of $\R^N$, with
\begin{equation}\label{subset3}
\Gamma_t\subset\bigcup\limits_{1\leq i\leq 
n}R_{i,t}\left(\{x=(x',x_N)\in\R^N:x'\in\omega_{i,t},\,x_N=\psi_{i,t}(x')\}\right).
\end{equation}
\begin{defi}[\!\!\cite{bh1,bh2}]\label{tran}
A transition front of problem \eqref{eq1.1} connecting $1$ and $0$ is a classical solution 
$u:\mathbb R\times\overline\Omega\to\mathbb R$  for which there 
exist some subsets $\{\Omega_t^\pm\}_{t\in\R}$ and $\{\Gamma_t\}_{t\in\R}$ 
satisfying \eqref{subset1}, \eqref{subset2} and \eqref{subset3} such that 
for any $\varepsilon>0$, there is a positive constant $M_\varepsilon>0$ satisfying
\begin{equation*}\label{tran1}
\left\{\begin{array}{ll}
\forall t\in\mathbb R,\ \forall x\in\Omega_t^+,~~(d_\Omega( x,\Gamma_t)
\geq M_\varepsilon)\Rightarrow (|u(t,x)- 1|\leq\varepsilon),\\
\forall t\in\mathbb R,\ \forall x\in\Omega_t^-,~~(d_\Omega(x,\Gamma_t)\geq 
M_\varepsilon) \Rightarrow (|u(t,x)- 0|\leq\varepsilon).
\end{array}\right.
\end{equation*}
Moreover, $u$ is said to have a global mean speed $\gamma(\geq0)$ if
\begin{equation*}\label{g-speed}
\frac{d_\Omega(\Gamma_t,\Gamma_s)}{|t-s|}\to\gamma\text{ as }|t-s|\to\infty.
\end{equation*}
\end{defi}

Propagation in unbounded domains depends on both the reaction and the
geometry of the domain. 
For bistable equations, changes in the geometry can lead to complete
invasion, partial invasion, or blocking. Berestycki, Hamel and Matano
\cite{bhm} studied planar fronts interacting with a compact obstacle and the
influence of the obstacle on their subsequent propagation.
In that homogeneous setting, their results distinguish recovery of the
propagating profile from uniform recovery throughout the exterior domain,
which depends on the obstacle geometry.
More recently,
Berestycki, Hamel and Matano \cite{bhm-wall} considered bistable fronts
encountering a perforated wall of finite thickness and infinite transverse
extent. They established a dichotomy between propagation through the wall
and blocking, and gave geometric sufficient conditions for each outcome.
In particular, sufficiently narrow holes cause blocking under transverse
periodicity or localization assumptions, whereas suitable large holes,
small-capacity walls, and parallel-blade walls permit propagation.
In cylinders with
varying cross section, Berestycki, Bouhours and Chapuisat
\cite{bbc-cylinders} gave geometric conditions ensuring complete invasion
and showed that an abrupt opening can block propagation. Ducasse and Rossi
\cite{dr-blocking} studied invasion and blocking in periodically perforated
domains, including examples where invasion occurs in one direction but is
blocked in the opposite direction. In domains consisting of a straight part
and a conical part, Hamel and Zhang \cite{hz-funnels} established a dichotomy
between blocking and spreading for bistable fronts, with sufficient
conditions involving the width of the straight part and the opening angle.
Guo, Hamel and Sheng
\cite{ghs-domains} proved that bistable transition fronts with complete
propagation in exterior domains have the same global mean speed as the
planar front. In domains with multiple cylindrical branches, they showed
that entire solutions emanating from some branches and propagating
completely are transition fronts with that same mean speed. 

Ignition reactions in exterior domains provide a further comparison
with the bistable setting considered here. Jia, Wang and Bu
\cite{jwb-combustion} proved, for time-periodic combustion equations, the
existence and uniqueness of an entire solution originating from a planar
front, complete invasion independently of the obstacle geometry within
their admissible class, and recovery of the same front after passage past
the obstacle. For spatially periodic combustion equations, Jia, Ma, Bao and
Guo \cite{jmbg-combustion} established analogous existence, uniqueness,
complete invasion and front recovery results for pulsating front-like
solutions. These results motivate the distinction between establishing
complete invasion and proving recovery of the incoming front.
Other related results can be referred to
\cite{bn-cylinders,cg-opening,guo-branches,heinze-perforated,ll-branches,
pauthier-cylinders,slww-obstacle,sww-cylinders,sw-systems,ys-combustion}.

We now recall  the
existence and uniqueness of the entire solution  established by Jia, Sheng and Wang
\cite[Theorem 1.2]{jsw}.

\begin{theo}\lb{existence}
Assume that (f1)-(f3)  as well as (A1)-(A2) hold.
Then there exists a unique entire solution \(u(x,t)\) of \eqref{eq1.1} in \(\Omega=\mathbb{R}^{N} 
\backslash K\),
such that \(0<u(x,t)<1\) and \(u_{t}(x,t)>0\) for all \((x,t)\in\bar\Omega\times\mathbb R\) and
\[
u(x,t)-U_1\left(x,x_1-c_1 t\right) \rightarrow 0\ \ \text{as  $t \rightarrow - \infty$   uniformly in $x \in \bar{\Omega}$.}
\]
\end{theo}

Our main result is Theorem \ref{convergence}, which establishes recovery
of the incoming front under complete propagation. We first state two
spatial asymptotic properties used. The first concerns
transverse translations escaping to infinity.

\begin{theo}\lb{intermediate}
Let $u(x,t)$ be the entire solution in Theorem \ref{existence}.  For any
sequence $(x_n')_{n\in\mathbb N}\subset\mathbb R^{N-1}$ with
$|x_n'|\to+\infty$, write $x_n'=\hat x_n'+\tilde x_n'$, where
$\hat x_n'\in[0,L_2)\times\cdots\times[0,L_N)$ and
$\tilde x_n'\in L_2\mathbb Z\times\cdots\times L_N\mathbb Z$. Then
\bes
 {u(x+(0,\tilde x_n'),t)\to U_1(x,x_1-c_1 t)
\quad\text{as }n\to+\infty}
\ees
locally uniformly in $(x,t)\in\R^N\times\R$.
\end{theo}

The second describes the behavior at spatial infinity for each fixed time.

\begin{theo}\lb{larget}
Let $u(x,t)$ be the entire solution of \eqref{eq1.1} given in Theorem \ref{existence}. 
Then 
  for every fixed $T\in\mathbb R$, there holds
\begin{equation}\label{t16-fixed-time-farfield}
 \lim_{R\to\infty}\ \sup_{\substack{x\in\overline\Omega,\ |x|\geq R}}
 |u(x,T)-U_1(x,x_1-c_1T)|=0.
\end{equation}
\end{theo}

We can now state the recovery result.

\begin{theo}\lb{convergence}
Let $u(x,t)$ be the entire solution of \eqref{eq1.1} given in Theorem \ref{existence}. 
Assume that the propagation is complete in the sense that 
\be\lb{com}
\lim_{t\to+\i} u( x,t)=1 \ \ \text{locally uniformly in $x\in\Omega$}.
\ee
Then there holds
\[
u(x,t)-U_1\left(x,x_1- c_1 t\right) \rightarrow 0\ \ \text{as  $t \rightarrow + \infty$   uniformly in
$x \in \bar{\Omega}$.}
\]
\end{theo}

\begin{rem}\label{rem1}
{
Theorem \ref{convergence} retains the complete propagation assumption
\eqref{com} but does not require \eqref{add}. The limiting front has the same
phase as the incoming front. In the homogeneous setting studied in
\cite{bhm}, geometric conditions such as star-shapedness or directional
convexity provide examples of complete invasion.
\par}
\end{rem}

The rest of this paper is organized as follows. In Section 2, we prove the intermediate behavior,
that is we prove Theorem \ref{intermediate}. In Section 3, we prove Theorem \ref{larget}. In Section 4, we 
prove Theorem \ref{convergence}.

\SE{Intermediate behavior}

In this section, we are concerned with the behavior near the horizon of a time global solution of 
\eqref{eq1.1},
which behaves like a pulsating traveling front for very negative time.
Here the horizon is referred to the limit that $|x'| \to+\i$  with $x'=(x_2,\cdots,x_N)$
when $x_1$ stays bounded. Let $(x_n')_{n\in\mathbb N}\subset \mathbb R^{N-1}$ 
be any sequence such that \(\left|x_{n}' \right| \rightarrow+\infty\) as \(n \rightarrow+\infty\). 
Then one can write 
$x_n'=\hat{x}_n'+\tilde x_n'$ with  
$\tilde x_n'\in L_2 \mathbb Z\times \cdots L_N \mathbb Z$. 

\vspace{0.3cm}

\begin{pr}[Proof of Theorem \ref{intermediate}]
Write  $W(x,t)=U_1(x,x_1-c_1t)$.
Let $x_n'=\hat x_n'+\tilde x_n'$ be the decomposition in the
statement. Since \( |\hat x_n'|\leq\left(\sum_{i=2}^N L_i^2\right)^{1/2}\), then 
\( |\tilde x_n'|\geq |x_n'|-|\hat x_n'|\longrightarrow+\infty\). Set
\[
 u_n(x,t)=u(x+(0,\tilde x_n'),t),
 \quad
 (x,t)\in(\Omega-(0,\tilde x_n'))\times\mathbb R.
\]

\medskip\noindent\textit{Step 1.}
Since \(0\leq u_n\leq1\), and the 
spatial periodicity gives
\[
 (u_n)_t=\Delta u_n+f(x,u_n)\ \  
 \text{in }(\Omega-(0,\tilde x_n'))\times\mathbb R.
\]
It then follows from the standard  parabolic 
and a diagonal
subsequence that $u_n\to v$ in $C^{2,1}_{\mathrm{loc}}(\mathbb{R}^N\times\mathbb{R})$ as \(n\to+\infty\), 
and \(0\leq v\leq1\) satisfies
\[
 v_t=\Delta v+f(x,v)\quad\text{in }\mathbb R^N\times\mathbb R.
\]
Write $e(t)=\sup_{\overline\Omega}|u(\cdot,t)-W(\cdot,t)|$.
Theorem \ref{existence} gives $e(t)\to0$ as $t\to-\infty$.
For each fixed $(x,t)$, 
there holds
\[
 |u_n(x,t)-W(x,t)|
 =|u(x+(0,\tilde x_n'),t)-W(x+(0, \tilde x_n'),t)|\leq e(t).
\]
Hence it follows that 
\begin{equation}\label{eq2.1}
 \|v(\cdot,t)-W(\cdot,t)\|_{L^\infty(\mathbb R^N)}
 \leq e(t)\longrightarrow0\quad \text{as \(t\to-\infty\).}
\end{equation}

\medskip\noindent\textit{Step 2.}
Let $L_f$ be a global Lipschitz constant of the extension of $f$. 
By \eqref{eq1.2},  
choose $0<\delta\leq\eta_0<1/4$ and $\omega>0$ so that
\[
 f_u(x,u)\leq-4\omega\ \ \text{for all \(u\leq2\delta\) or \(u\geq1-2\delta.\)}
\]
Take a constant \(C>0\) such that 
\[
U_1(x,s)\leq\delta \ \ \text{for all \(s\geq C\) \ \ and}\ \  U_1(x,s)\geq 1-\delta \ \ \text{for all \(s\leq -C\) }
.\]
It follows from \cite[Lemma 2.5]{gllw} that there is a \(\kappa>0\) such that
\[
 \kappa:=\inf_{\delta\leq U_1(x,s)\leq1-\delta}(-\partial_s U_1(x,s)).
\]
Choose $\rho\geq(L_f+\omega)/(\kappa c_1)$. Fix $0<\varepsilon\leq\delta$.
by \eqref{eq2.1}, there is $T_\varepsilon$ so that 
\[
 \|v(\cdot,t)-W(\cdot,t)\|_{L^\infty(\mathbb R^N)}
 \leq \varepsilon\ \ \text{ for all $t\leq T_\varepsilon$.}
\]
For an arbitrary $t_0\leq T_\varepsilon$, define
\[
 u^\pm(x,t)=U_1\left(x,x_1-c_1(t\pm \frac{\rho\varepsilon}{\omega}(1-e^{-\omega(t-t_0)}))\right)\pm \varepsilon e^{-\omega(t-t_0)}\ \ 
\text{for $t\geq t_0$ and \(x\in\mathbb{R}^N\).}
\]
Let $\mathcal N[  u^\pm]= u^\pm_t-\Delta  u^\pm-f(x, u^\pm)$.
After some direct calculations, one has 
\begin{align*}
 \mathcal N[u^-]&=c\rho \varepsilon e^{-\omega(t-t_0)} \partial_\xi U_1+\omega \varepsilon e^{-\omega(t-t_0)}+f(x, U_1)-f(x,u^-).
\end{align*}
If  \(s\geq C\) or \(s\leq -C\), then 
$U_1\leq\rho$ or $U_1\geq1-\rho$, then $u^- \leq \rho$ or $u^- \geq1-2\rho$. Hence
\[
 f(x,U_1)-f(x,u^-) \leq-4\omega  \varepsilon e^{-\omega(t-t_0)},
\]
whence
\[
\mathcal N[u^-]\leq-3\omega  \varepsilon e^{-\omega(t-t_0)}\leq0.
\]
If \(-C\leq s\leq C\), then 
\(f(x,U_1)-f(x,u^-) \leq L_f  \varepsilon e^{-\omega(t-t_0)},\)
hence 
\[ 
\mathcal N[u^-]\leq -c\rho \kappa \varepsilon e^{-\omega(t-t_0)} +\omega \varepsilon e^{-\omega(t-t_0)}
+L_f  \varepsilon e^{-\omega(t-t_0)}\leq0.
\]
At $t=t_0$,   it holds
\[
 u^-(x,t_0)=U_1(x,x_1-c_1 t_0)-\varepsilon\leq v(x,t_0)
 \leq U_1(x,x_1-c_1t_0)+\varepsilon=u^+(x,t_0).
\]
It then follows from the comparison principle that 
\begin{equation}\label{eq2.2}
u^-(x,t)\leq v(x,t)\leq u^+(x,t) \ \ \text{for \(x\in\mathbb{R}^N\) and \(t\geq t_0\).}
\end{equation}
 Let $t_0\to-\infty$ in
\eqref{eq2.2}. 
Then 
\[
 U_1(x,x_1-c_1 t+c_1\rho\varepsilon/\omega)\leq v(x,t)
 \leq U_1(x,x_1-c_1 t-c_1\rho\varepsilon/\omega)  \ \ \text{for \(x\in\mathbb{R}^N\) and \(t\in \R\).}
\]
Notice that the constants $\rho,\omega$ are independent of $\varepsilon$.
Letting $\varepsilon\to0$ gives $v(x,t)=U(x,x_1-c_1t)$.

\medskip\noindent\textit{Step 3.}
Let us finally show that, for every compact set
$Q\subset\mathbb R^N\times\mathbb R$ and every $\varepsilon>0$, there
exists $n(Q,\varepsilon)$ such that, for all $n\geq n(Q,\varepsilon)$,
\[
 Q\subset(\Omega-(0, \tilde x_n'))\times\mathbb R
 \ \ \text{and}\ \  \sup_{(x,t)\in Q}|u_n(x,t)-W(x,t)|<\varepsilon.
\]
If it is not true, then there exist a fixed
compact $Q$, a number $\varepsilon_0>0$, and a strictly increasing sequence
of indices $n_j\to\infty$ such that 
\begin{equation*}
 \sup_{(x,t)\in Q}|u_{n_j}(x,t)-W(x,t)|\geq\varepsilon_0
 \quad\text{for every }j.
\end{equation*}
It follows from the proof of Step 1 that there 
are a subsequence \( u_{n_{j_\ell}}\) and  an entire solution $v_*$ such that
\[
 u_{n_{j_\ell}}\longrightarrow v_*
 \quad\text{in }C^{2,1}_{\mathrm{loc}}(\mathbb R^N\times\mathbb R)
\]
and
\[\|v_*(x,t)-U_1(\cdot, x_1-t)\|_\infty\leq e(t)\ \ \text{for every $t$.}\] 
In addition, it can be inferred from Step 2 that 
\[
 U_1(x,x_1-c_1t +c_1 \rho\varepsilon/\omega)\leq v_*(x,t)
 \leq U_1(x,x_1-c_1t -c_1 \rho\varepsilon/\omega)
\]
for every $\varepsilon\in(0,\delta]$ and every $(x,t)$,
Letting $\varepsilon\downarrow0$ yields $v_*(x,t)=U_1(x,x_1-c_1t)$.
Local uniform convergence now gives
\[
 \begin{aligned}
 \varepsilon_0
 &\leq\sup_{(x,t)\in Q}|u_{n_{j_\ell}}(x,t)-W(x,t)|\\
 &\leq\sup_{(x,t)\in Q}|u_{n_{j_\ell}}(x,t)-v_*(x,t)|
       +\sup_{(x,t)\in Q}|v_*(x,t)-W(x,t)|\\
 &=\sup_{(x,t)\in Q}|u_{n_{j_\ell}}(x,t)-v_*(x,t)|
 \longrightarrow0.
 \end{aligned}
\]
This contradicts $\varepsilon_0>0$. Hence the full sequence $u_n$ converges
to $U$ locally uniformly on $\mathbb R^N\times\mathbb R$. This completes the proof.
\end{pr}

\SE{Large time behavior}

In this section, we prove Theorem \ref{larget}: at each fixed time,
the entire solution approaches the incoming pulsating front uniformly
outside sufficiently large balls. We first control the solution in the
two longitudinal tails and then apply Theorem \ref{intermediate} to the
remaining slab as the transverse variable tends to infinity. Combining
these estimates yields the asserted spatial far field convergence.
\vspace*{0.3cm}

\begin{pr}[Proof of Theorem \ref{larget}]
\medskip\noindent\textit{Step 1.}   
Let $W$ and $L_f$ be given in the proof of Theorem \ref{intermediate}.
Fix an arbitrary finite time $T\in\mathbb R$ and an arbitrary
$\varepsilon>0$. Choose
\[
 0<\eta<\min\{\varepsilon/4,\eta_0,1/4\}\ \ \text{and}\ \  \lambda>0.
\]
By (f3) and $f(x,0)=0$, one has \(f(x,\eta)<0\).
In particular, $|f(x,r)-f(x,s)|\leq L_f|r-s|$  for all $x\in\Omega$ and $r,s\in\mathbb R$.
The past asymptotics in Theorem \ref{existence} implies that there is
$\tau\leq T$ such that
\begin{equation}\label{eq3.1}
 \sup_{x\in\overline\Omega}|u(x,\tau)-W(x,\tau)|\leq\eta/2.
\end{equation}
Choose $C>0$ such that
\[
 U_1(x,\xi)\geq1-\eta/2\ \ \text{for \(\xi\leq-C\)\ \ and}\ \ 
 U_1(x,\xi)\leq\eta/2\ \ \text{for \(\xi\geq C\)}
\]
uniformly in $x$. Choose $R_K>0$ such that $K\subset B_{R_K}$. Let 
\[
 a_-<\min\{-R_K, c_1\tau-C\},\ \ 
 b>\max\{R_K, c_1\tau+C\}.
\]
Then $K\subset\{a_-<x_1<b\}\times \mathbb{R}^{N-1}$ and
\[
 W(x,\tau)\geq1-\eta/2\ \ \text{for \(x_1\leq a_-\) \ \ and }\ \ 
 W(x,\tau)\leq\eta/2\ \ \text{for \(x_1\geq b\)}
\]
uniformly in \(x'\in\mathbb{R}^{N-1}\).

\medskip\noindent\textit{Step 2.}
For $x_1\leq a_-$, by \eqref{eq3.1} one has that 
\[
 u(x,\tau)\geq W(x,\tau)-\eta/2\geq1-\eta.
\]
Since $u_t>0$, $W_t>0$ and $T\geq\tau$, then
\[
 1-\eta\leq u(x,\tau)\leq u(x,T)\leq1,
 \ \ 
 1-\eta/2\leq W(x,\tau)\leq W(x,T)\leq1.
\]
As a result, there holds
\begin{equation}\label{eq3.2}
 |u(x,T)-W(x,T)|\leq\eta<\varepsilon\ \ 
 \text{for \(x_1\leq a_-\).}
\end{equation}

\medskip\noindent\textit{Step 3.}
For $x_1\geq b$ and $t\geq\tau$, define
\[
u^+(x,t)=\eta+\exp\{-\lambda(x_1-b)+(\lambda^2+L_f)(t-\tau)\}.
\]
After some calculations, one has
\[\begin{aligned}
u^+_t-\Delta u^+-f(x,u^+)
= L_f \exp\{-\lambda(x_1-b)+(\lambda^2+L_f)(t-\tau)\}-f(x,u^+).
\end{aligned}
\]
Since 
\[
 f(x,u^+)\leq f(x,\eta)+L_f(u^+-\eta)\leq L_f e^{-\lambda(x_1-b)+(\lambda^2+L_f)(t-\tau)},
\]
then 
\[
u^+_t-\Delta u^+-f(x,u^+)\geq 0,\ \ \text{\(x_1>b\).}
\]
The extension of $f$ is used here because $u^+$ may exceed $1$. 
At the initial time, one has  from \eqref{eq3.1} that
\[
 u(x,\tau)\leq W(x,\tau)+\eta/2\leq\eta
 \leq {u^+(x,\tau)}\ \ \text{for \( x_1\geq b\).}
\]
On the boundary $x_1=b$, there holds
\[
u^+(x,t)=\eta+e^{(\lambda^2+L_f)(t-\tau)}\geq1\geq u(x,t) 
 \ \ \text{for \(t\geq\tau\).}
\]
It then follows from the comparison principle that 
\begin{equation*}
 0\leq u(x,t)\leq\eta+
 e^{-\lambda(x_1-b)+(\lambda^2+L_f)(t-\tau)}
 \ \ \text{for \(x_1\geq b\) \ and\  \(t\geq\tau\).}
\end{equation*}
Choose $a_+\geq b$ so that
\[
 a_+\geq b+\frac{\lambda^2+L_f}{\lambda}(T-\tau)
                  +\frac1\lambda\log(1/\eta)
 \ \ \text{and}\  \  a_+\geq {c_1}T+C.
\]
For $x_1\geq a_+$, one obtains that
\[
 -\lambda(x_1-b)+(\lambda^2+L_f)(T-\tau)\leq\log\eta,
 \quad W(x,T)\leq\eta/2,\quad u(x,T)\leq2\eta.
\]
Since both functions are nonnegative, then one gets
\begin{equation}\label{eq3.3}
 |u(x,T)-W(x,T)|\leq\max\{u(x,T),W(x,T)\}\leq2\eta<\varepsilon
 \ \ \text{for \(x_1\geq a_+\).}
\end{equation}

\medskip\noindent\textit{Step 4.}
We prove that there exists $B>R_K$ such that
\begin{equation}\label{eq3.4}
 |u(x_1,y,T)-W(x_1,y,T)|\leq\varepsilon
 \ \ \text{ for \(a_-\leq x_1\leq a_+,\ |y|\geq B\).}
\end{equation}
If this is not true, then, for each integer $j\geq1$, there are
$x_{j,1}\in[a_-,a_+]$ and $y_j\in\mathbb R^{N-1}$ such that \(|y_j|\geq R_K+j\)
and \(|u(x_{j,1},y_j,T)-W(x_{j,1},y_j,T)|\geq  {\varepsilon}.\)
Decompose $y_j=k_j'+\widehat y_j$, where
\(
 k_j'\in L_2\mathbb Z\times\cdots\times L_N\mathbb Z\) and \(\widehat y_j\in\prod_{i=2}^N[0,L_i). \)
Then \( |k_j'|\geq|y_j|-\left(\sum_{i=2}^N L_i^2\right)^{1/2}\longrightarrow\infty\) as \(j\to\infty\).
The points $(x_{j,1},\widehat y_j)$ all belong to
$[a_-,a_+]\times\prod_{i=2}^N[0,L_i]$.
Theorem \ref{intermediate} therefore gives
\[
 \sup_{\substack{x_1\in[a_-,a_+]\\y\in\prod_{i=2}^N[0,L_i]}}
 |u((x_1,y)+(0,k_j'),T)-W((x_1,y),T)|\longrightarrow0 \ \ \text{as \(j\to+\infty\).}
\]
The periodicity further implies that  $W(x_{j,1},y_j,T)=W(x_{j,1},\widehat y_j,T)$. 
As a result, one derives that 
\[
 \begin{aligned}
  {\varepsilon} 
 &\leq |u(x_{j,1},y_j,T)-W(x_{j,1},y_j,T)|\\
 &=|u((x_{j,1},\widehat y_j)+(0,k_j'),T)
    -W(x_{j,1},\widehat y_j,T)|\\
 &\leq\sup_{\substack{x_1\in[a_-,a_+]\\y\in\prod_{i=2}^N[0,L_i]}}
 |u((x_1,y)+(0,k_j'),T)-W((x_1,y),T)|\longrightarrow0\ \ \text{as \(j\to+\infty\),}
 \end{aligned}
\]
which is impossible. This proves \eqref{eq3.4}.

Define
\[
 C_{\varepsilon,T}
 =\bigl([a_-,a_+]\times\overline B_B^{N-1}\bigr)\cap\overline\Omega.
\]
If $x=(x_1,y)\in\overline\Omega\setminus C_{\varepsilon,T}$, then either
$x_1<a_-$, or $x_1>a_+$, or $a_-\leq x_1\leq a_+$ and $|y|>B$.
It then can be deduces from estimates \eqref{eq3.2}, \eqref{eq3.3}, and \eqref{eq3.4}
that
\[
 |u(x,T)-W(x,T)|\leq\varepsilon
 \ \ {x\in\overline\Omega\setminus C_{\varepsilon,T}}.
\]
Since $T$ and $\varepsilon$ were arbitrary, this yields
\eqref{t16-fixed-time-farfield}. The proof is complete.
\end{pr}

\vspace*{0.3cm}

\SE{Recovery of the pulsating front after complete propagation}

In this section, we prove Theorem \ref{convergence}, which asserts
uniform convergence to the incoming pulsating front under the assumption
of complete propagation \eqref{com}. We first obtain bounds in terms of time shifts
of the front and then use a positive adjoint weight to prove decay of
localized phase perturbations. Finally, comparison with suitable super-sub solutions, together with Theorem \ref{larget}, yields the
desired convergence.

 Let
\( x=(x_1,y)\)  and \(a_+=\max\{a,0\}.\)
Write 
\[
 W(x,t)=U_1(x,x_1-c_1t),\quad x\in\mathbb R^N,\ t\in\mathbb R.
\]
For each $r\in\mathbb R$, define
\[
 W^r(x,t)=W(x,t+r)=U_1(x,x_1-c_1t-c_1r),\quad x\in\mathbb R^N,\ t\in\mathbb R.
\]
The superscript $r$ is a time shift, not a power.  In the sequel, we use
the letter $C$ denotes a positive constant that may change between inequalities. 
Extend $f$ to all $x\in\mathbb R^N$ and $u\in\mathbb R$ by  \eqref{eq1.6}.
Let  $L_f$ be given as in the proof of Theorem \ref{intermediate}.
From (f3) we choose $\omega$ and \(\eta_0\) so that
\begin{equation}\label{eq4.1}
 \begin{gathered}
 0< \eta_0<1/4,\ \  0<\omega\leq\beta_0/8,\\
 |f(x,v)-f(x,w)|\leq L_f|v-w|,\ \ 
 f_u(x,v)\leq-4\omega,\ \   {v\leq2 \eta_0\text{ or }v\geq1-2 \eta_0}.
 \end{gathered}
\end{equation}
Let $P=L_1/c_1$. 
Notice that 
\begin{equation*}
 W_t=-c_1\partial_sU_1(x,s) >0.
\end{equation*}
It is obvious that there is $ C_1>0$  so that
\[
  \eta_0\leq U_1(x,s)\leq1- \eta_0 \quad \text{for \(|s|\leq  {C_1}\).}
\]
Define
\begin{equation}\label{eq4.2}
  {\kappa}:=\inf_{ \eta_0\leq U_1(x,s)\leq1- \eta_0}(-c_1\partial_sU_1(x,s))>0,
 \ \  M:=\sup_{(x,t)\in\mathbb{R}^N\times\mathbb{R}}W_t(x,t)<\infty,
 \ \   \rho:=\frac{L_f+\omega}{ \kappa}.
\end{equation}
Fix $R>0$ with $K\subset B_R$. Choose $ \mu>0$ such that
$ {\mu}^2+ \mu(N-1)/R\leq\omega$. Write \( \mathcal D=\mathbb R^N\setminus\overline B_R\).
Define
\begin{equation}\label{eq4.3}
 \zeta(x)=e^{- \mu(|x|-R)}, \ \   x\in\overline{\mathcal D}.
\end{equation}
Then \(0<\zeta\leq1\), \(\zeta|_{\partial B_R}=1\), \(\Delta\zeta=\left( {\mu}^2-\frac{ \mu(N-1)}{|x|}\right)\zeta\)
and \(|\Delta\zeta|\leq\omega\zeta\).
In addition, there holds
\begin{equation}\label{eq4.4}
 \zeta(x)\leq\min\{1,e^{ \mu(R+ C_1)- \mu c_1(t+r)}\}
 \quad\text{when } {\eta_0}\leq W^r(x,t)\leq1- \eta_0.
\end{equation}
Set
\[
 \Nop[v]=v_t-\Delta v-f(x,v).
\]
For $v\in[0,1]$ and $0\leq r\leq {\eta_0}$, if $v\leq {\eta_0}$ or $v\geq1- \eta_0$,  by \eqref{eq4.1} one obtains that
\begin{equation}\label{eq4.5}
 \begin{aligned}
 f(x,v)-f(x,v-r)&=\int_{v-r}^v f_u(x,q)\,dq\leq-4\omega r,\\
 f(x,v)-f(x,v+r)&=-\int_v^{v+r} f_u(x,q)\,dq\geq4\omega r.
 \end{aligned}
\end{equation}

\begin{lem}\label{lem-t17-trapping}
Assume the hypotheses of Theorem \ref{convergence}. There exists a
constant $D\geq1$ such that
\begin{equation}\label{eq4.6}
 \|(W^{-D}(\cdot,t)-u(\cdot,t))_+\|_{L^\infty(\overline\Omega)}
 +\|(u(\cdot,t)-W^D(\cdot,t))_+\|_{L^\infty(\overline\Omega)}
 \longrightarrow0 \ \   \text{as}\  t\to\infty.
\end{equation}
\end{lem}
\begin{pr}
By the complete propagation \eqref{com}, one obtains that
\begin{equation}\label{eq4.7}
 \sup_{x\in \mathcal C}|1-u(x,t)|\longrightarrow0
 \quad\text{for every compact set } \mathcal C\subset\overline\Omega.
\end{equation}
Choose
\[
 0<\delta\leq\min\{ \eta_0/2,1/(4 \rho)\}.
\]
By \eqref{eq4.7}, there exists a $T$ such that
\begin{equation}\label{eq4.8}
 u(x,t),W(x,t)\geq1-\delta, \ \   x\in\partial B_R,\ t\geq T.
\end{equation}
By Theorem \ref{larget}, there is a compact set $\mathcal C\subset\overline\Omega$
with $|u(x,T)-W(x,T)|\leq\delta$ outside $\mathcal C$.
On $\mathcal C\cap\overline{\mathcal D}$, $W^{-d_0}(x,T)\to0$ and
$W^{d_0}(x,T)\to1$ uniformly as $d_0\to\infty$. Choose $d_0\geq1$
so that $W^{-d_0}\leq\delta$ and $W^{d_0}\geq1-\delta$ there.
Note that $W^{-d_0}\leq W\leq W^{d_0}$. Hence
\begin{equation}\label{eq4.9}
 W^{-d_0}(x,T)-\delta\leq u(x,T)\leq W^{d_0}(x,T)+\delta
 , \ \   x\in\overline{\mathcal D}.
\end{equation}
For $t\geq T$,  define 
\begin{equation}\label{eq4.10}
 d'(t)= \rho\delta\left[e^{-\omega(t-T)}
       +\min\{1,e^{ \mu(R+ C_1)- \mu c_1(t-d)}\}\right],
 \ \  d(T)=d_0.
\end{equation}
one has that \( 0\leq d'\leq2 \rho\delta\leq\tfrac12\) and 
\(t-d(t)\geq T-d_0+\tfrac12(t-T).\)
Consequently,
\[
 \int_T^\infty\min\{1,e^{ \mu(R+ C_1)- \mu c_1(t-d(t))}\}\,dt
 \leq\frac{2}{ \mu c_1}e^{ \mu(R+ C_1)- \mu c_1(T-d_0)}<\infty.
\]
Integrating \eqref{eq4.10} gives
\begin{equation*}
 D:=\lim_{t\to\infty}d(t)
 \leq d_0+\frac{ \rho\delta}{\omega}
 +\frac{2 \rho\delta}{ \mu c_1}e^{ \mu(R+ C_1)- \mu c_1(T-d_0)}<\infty.
\end{equation*}
On $\overline{\mathcal D}\times[T,\infty)$, define
\[
 w^\mp(x,t)=W(x,t\mp d(t))\mp \delta e^{-\omega(t-T)}\mp \zeta(x).
\]
After some computations, one gets that 
\[
 \begin{aligned}
 \Nop[w^-]=-d'(t)W_t(x,t-d(t))
       +\omega\delta e^{-\omega(t-T)}+\delta\Delta\zeta(x)+f(x,W(x,t-d(t)))-f(x,w^-(x,t)).
 \end{aligned}
\]
Since $|\Delta\zeta|\leq\omega\zeta$ and
$0\leq\delta e^{-\omega(t-T)}+\delta\zeta\leq2\delta\leq {\eta_0}$,
if
$W(x,t-d(t))\leq {\eta_0}$ or $W(x,t-d(t))\geq1- \eta_0$, it follows from \eqref{eq4.5} that
\[
 \Nop[w^-]\leq-d'W_t(x,t-d)
       -3\omega\delta\bigl(e^{-\omega(t-T)}+\zeta\bigr)\leq0.
\]
If $ \eta_0\leq W(x,t-d(t))\leq1- \eta_0$, then
\eqref{eq4.4} and \eqref{eq4.2} give
\[
 \begin{aligned}
 \Nop[w^-]
 &\leq-d' \kappa+(L_f+\omega)\delta
                       \bigl(e^{-\omega(t-T)}+\zeta\bigr)\\
 &\leq(- \rho {\kappa}+L_f+\omega)\delta
 \left(e^{-\omega(t-T)}
       +\min\{1,e^{ \mu(R+ C_1)- \mu c_1(t-d)}\}\right)=0.
 \end{aligned}
\]
On $\partial B_R$, it can be deduced from \eqref{eq4.8} and $W\leq1$ that
\[
 w^-(x,t)=W(x,t-d(t))-\delta e^{-\omega(t-T)}-\delta
 \leq1-\delta\leq u(x,t).
\]
At $t=T$, it infers from \eqref{eq4.9} that
\[
 w^-(x,T)=W(x,T-d_0)-\delta-\delta e^{- \mu(|x|-R)}
 \leq W^{-d_0}(x,T)-\delta\leq u(x,T).
\]
It then follows from the comparison principle that
\begin{equation*}
 W^{-d(t)}-\delta e^{-\omega(t-T)}-\delta e^{- \mu(|x|-R)}
 \leq u
 \leq W^{d(t)}+\delta e^{-\omega(t-T)}+\delta e^{- \mu(|x|-R)}
\end{equation*}
in $\overline{\mathcal D}\times[T,\infty)$.
Since $W^{-D}\leq W^{-d(t)}$ and $W^{d(t)}\leq W^D$, for $H>R$, one arrives at
\[
 \sup_{|x|\geq H}\max\{(W^{-D}-u)_+,(u-W^D)_+\}
 \leq\delta e^{-\omega(t-T)}+\delta e^{- \mu(H-R)}.
\]
On $\overline\Omega\cap\overline B_H$, \eqref{eq4.7} gives
$u\to1$ uniformly. Then \(x_1-c_1t\mp c_1D\leq H-c_1t+c_1D\longrightarrow-\infty\)
as \(t\rightarrow+\infty\), hence \(W^{\pm D}\longrightarrow1\) uniformly on $\overline\Omega\cap\overline B_H$.
Thus, for every $H>R$,
\[
 \limsup_{t\to\infty}\sup_{x\in\overline\Omega}
 \max\{(W^{-D}-u)_+,(u-W^D)_+\}
 \leq\delta e^{- \mu(H-R)}.
\]
Letting $H\to\infty$ proves \eqref{eq4.6}. This completes the proof.
\end{pr}
\vspace{0.3cm}

Set
\[
 X=(s,y),\quad {s}=x_1-c_1t,\quad
 \mathbb T'=\prod_{i=2}^N(\mathbb R/L_i\mathbb Z),\quad
 \mathcal C=\mathbb R\times\mathbb T',
 \quad e_{s}=(1,0,\ldots,0).
\]
For each constant $r\in\mathbb R$, define
\begin{equation}\label{eq4.11}
 V^r(X,t)=U_1((s+c_1t,y),s-c_1r), \ \ 
 p^r=\partial_rV^r>0, \ \  b^r(X,t)=f_u((s+c_1t,y),V^r(X,t)).
\end{equation}
For a function $\psi^0$ on $\mathcal C\times\mathbb R$, write
\[
 \psi^r(s,y,t)=\psi^0(s-c_1r,y,t+r), \ \  r\in\mathbb R.
\]

\begin{lem}\label{lem-t17-weight}
There exists $\psi^0>0$ defined on $\mathcal C\times\mathbb R$, 
periodic in $y$ and $P$-periodic in $t$, and 
\( \psi^r(s,y,t)=\psi^0(s-c_1r,y,t+r)\) with \(r\in\mathbb R\)
satisfies 
\begin{equation}\label{eq4.12}
 -\psi_t^r=\Delta_X\psi^r-c_1\psi_{s}^r+b^r\psi^r\ \ \text{and}
 \ \  \int_{\mathcal C}p^r(X,t)\psi^r(X,t)\,dX=1.
\end{equation}
\end{lem}

\begin{pr}
By (f1)--(f2),  there holds
\begin{equation*}
 |f_u(x,v)-f_u(\widetilde x,\widetilde v)|
 \leq C|x-\widetilde x|^{\alpha/2}+C|v-\widetilde v|.
\end{equation*}
For $|X-\widetilde X|$, $|t-\widetilde t|\leq1$, one has
\[
 \begin{aligned}
 |b^0(X,t)-b^0(\widetilde X,\widetilde t)|
 \leq C\bigl(|X-\widetilde X|^{\alpha/2}
  +|t-\widetilde t|^{\alpha/4}\bigr).
 \end{aligned}
\]
 Also $|b^0(X,t)|\leq L_f$ for every $(X,t)$.
By \eqref{eq4.11}, there hold
\[
 \begin{aligned}
 &V^r(s,y,t)=W((s+c_1t,y),t+r),&\\
 &V_t^r=W_t((s+c_1t,y),t+r)
+c_1W_{x_1}((s+c_1t,y),t+r),&\\
 &V_{s}^r=W_{x_1}((s+c_1t,y),t+r), \ \ 
 \Delta_XV^r=\Delta_xW((s+c_1t,y),t+r),&\\
 &p^r(s,y,t)=W_t((s+c_1t,y),t+r)
 =-c_1\partial_sU_1((s+c_1t,y),s-c_1r)>0.&
 \end{aligned}
\]
Substitution into $W_t=\Delta_xW+f(x,W)$ gives
\begin{equation}\label{eq4.13}
 V_t^r=\Delta_XV^r+c_1V_{s}^r+f((s+c_1t,y),V^r), \ \ 
 p_t^r=\Delta_Xp^r+c_1p_{s}^r+b^rp^r.
\end{equation}
In addition,  there hold $0<p^r\leq M$.
Since $c_1P=L_1$ and $U_1$ is periodic, then
\[
 V^r(s,y,t+P)
 =U_1((s+c_1t+L_1,y),s-c_1r)=V^r(s,y,t).
\]
Also it holds 
\begin{equation}\label{eq4.14}
 \lim_{H\to\infty}\sup_{\substack{|s|\geq H,\ y\in\mathbb T',\ 0\leq t\leq P}}
 p^0(s,y,t)=0.
\end{equation}
Choose $Z_0>0$ so large that, for every $y,t$,
\[
V^0(s,y,t)\leq\eta_0\ \ \text{for \(s\geq Z_0\) and}\ \ 
V^0(s,y,t)\geq1-\eta_0\ \ \text{for \(s\leq-Z_0\).}
\]
Then \eqref{eq4.1} implies
\begin{equation}\label{eq4.15}
 b^0(X,t)\leq-4\omega, \ \   |s|\geq Z_0.
\end{equation}
Let $\mathcal X=C_0(\mathcal C)$, the real Banach space of continuous
functions vanishing uniformly as $|s|\to\infty$, with
\[
 \|g\|_\infty=\sup_{X\in\mathcal C}|g(X)|, \ \ 
 \|T\|=\sup_{\|g\|_\infty\leq1}\|Tg\|_\infty.
\]
For $Y=(\xi,\upsilon)$ and $\tau>0$, the heat kernels are
\[
 \begin{aligned}
 H_{\mathbb T'}(\tau;y,\upsilon)
 &=\prod_{i=2}^N\sum_{k\in\mathbb Z}
   \frac{\exp(-(y_i-\upsilon_i+kL_i)^2/(4\tau))}{\sqrt{4\pi\tau}},\\
 H_{c_1}(\tau;X,Y)
 &=\frac{\exp(-(s+c_1\tau-\xi)^2/(4\tau))}{\sqrt{4\pi\tau}}
 H_{\mathbb T'}(\tau;y,\upsilon).
 \end{aligned}
\]
It is obvious that 
\[
 \begin{gathered}
 H_{c_1}>0, \ \ 
 \int_{\mathcal C}H_{c_1}(\tau;X,Y)\,dY
 =\int_{\mathcal C}H_{c_1}(\tau;X,Y)\,dX=1,\\
 \int_{\mathcal C}H_{c_1}(a;X,Z)H_{c_1}(b;Z,Y)\,dZ
 =H_{c_1}(a+b;X,Y), \ \   a,b>0.
 \end{gathered}
\]
For initial time $\tau<t$ and integers $n\geq0$, define
\[
 \begin{aligned}
 \Gamma_0(t,\tau;X,Y)&=H_{c_1}(t-\tau;X,Y),\\
 \Gamma_{n+1}(t,\tau;X,Y)
 &=\int_{\tau}^t\!\int_{\mathcal C}H_{c_1}(t-a;X,Z)
(b^0(Z,a)+L_f)\Gamma_n(a,\tau;Z,Y)\,dZ\,da.
 \end{aligned}
\]
Since $0\leq b^0+L_f\leq2L_f$, induction using the convolution identity
gives
\[
 0\leq\Gamma_n(t,\tau;X,Y)
 \leq\frac{(2L_f)^n(t-\tau)^n}{n!}H_{c_1}(t-\tau;X,Y).
\]
As  a result,  the series
\[
 \Gamma(t,\tau;X,Y)=e^{-L_f(t-\tau)}\sum_{n=0}^{\infty}\Gamma_n(t,\tau;X,Y)
\]
defines a  continuous kernel.  Substitution into the integral equation shows that
\[
 E(t,\tau)g(X):=\int_{\mathcal C}\Gamma(t,\tau;X,Y)g(Y)\,dY
\]
solves the following equation
\[
 \begin{aligned}
 E(t,\tau)g(X)
 &=\int_{\mathcal C}H_{c_1}(t-\tau;X,Y)g(Y)\,dY\\
 &\quad+\int_{\tau}^t\!\int_{\mathcal C}H_{c_1}(t-a;X,Y)
     b^0(Y,a)E(a,\tau)g(Y)\,dY\,da.
 \end{aligned}
\]
It follows from the interior estimates \cite[Theorems 7.22 and 4.9]{lieberman-parabolic}
and Sobolev embedding  that 
\[
 \partial_t E(t,\tau)g=(\Delta_X+c_1\partial_{s}+b^0(X,t))E(t,\tau)g.
\]
The nonnegative series and its factorial bound give directly
\begin{equation*}
 0<e^{-L_f(t-\tau)}H_{c_1}(t-\tau;X,Y)
 \leq\Gamma(t,\tau;X,Y)
 \leq e^{L_f(t-\tau)}H_{c_1}(t-\tau;X,Y).
\end{equation*}
Hence
\begin{equation}\label{eq4.16}
 E(t,\tau)\geq0, \ \ 
 |E(t,\tau)g|\leq E(t,\tau)|g|, \ \ 
 \|E(t,\tau)\|\leq e^{L_f(t-\tau)}.
\end{equation}
For $g\in\mathcal X$, choose $L$ so that $|g|\leq\varepsilon$ on
$|\xi|>L$. Let $t>\tau$ be fixed, then one has
\[
 |E(t,\tau)g(X)|
 \leq e^{L_f(t-\tau)}\left[
 \frac{2L\|g\|_\infty}{\sqrt{4\pi(t-\tau)}}
 e^{-(|s|-L-c_1(t-\tau))_+^2/(4(t-\tau))}+\varepsilon\right],\]
hence
\[
 \limsup_{|s|\to\infty}\sup_{y\in\mathbb T'}|E(t,\tau)g(s,y)|
 \leq e^{L_f(t-\tau)}\varepsilon\rightarrow0\ \ \text{as \(\varepsilon\rightarrow0\)}.
\]
Thus $E(t,\tau)\mathcal X\subset\mathcal X$. 
For $h>0$,
\[
 \left\|E(\tau+h,\tau)g-
\int_{\mathcal C}H_{c_1}(h;\cdot,Y)g(Y)\,dY\right\|_\infty
 \leq(e^{L_fh}-1)\|g\|_\infty.
\]
Hence  $E(\tau+h,\tau)g\to g$ in $\mathcal X$ as $h\downarrow0$.
In addition, there hold
\[
 E(t,\tau)E(\tau,a)=E(t,a),\ \  E(\tau,\tau)=I, \ \ 
 E(t+P,\tau+P)=E(t,\tau),
\]
where $I$ is the identity operator on $\mathcal X$.  

Define
\[
 \widehat b=\min\{b^0,-2\omega\}, \ \  d_b=b^0-\widehat b\geq0,
 \ \  \mathcal M=E(P,0).
\]
Here $d_b(\cdot,t)$ also denotes multiplication by this function, and
$\|d_b\|_\infty=\sup_{X\in\mathcal C,\,0\leq t\leq P}|d_b(X,t)|$.
By \eqref{eq4.15}, $d_b(X,t)=0$ when $|s|\geq Z_0$.
Define $\widehat E(t,\tau)g$ by the same initial value problem with
$b^0$ replaced by $\widehat b$. Since $\widehat b\leq-2\omega$,
comparison with both $\pm e^{-2\omega(t-\tau)}\|g\|_\infty$ gives
\begin{equation}\label{eq4.17}
 \|\widehat E(t,\tau)\|\leq e^{-2\omega(t-\tau)}.
\end{equation}
Then for every $g\in\mathcal X$, one has 
\[
 \begin{gathered}
 (\partial_t-\Delta_X-c_1\partial_{s}-b^0(X,t))
    \{E(t,0)g-\widehat E(t,0)g\}=d_b(X,t)\widehat E(t,0)g,\\
 E(0,0)g-\widehat E(0,0)g=0.
 \end{gathered}
\]
Hence 
\[
 E(t,0)g-\widehat E(t,0)g
 =\int_0^t E(t,\tau)\{d_b(\cdot,\tau)\widehat E(\tau,0)g\}\,d\tau\ \ \text{for $0\leq {\tau}\leq t\leq P$ }.
\]
In particular, at $t=P$ we obtain
\begin{equation*}
 \mathcal M-\widehat E(P,0)
 =\int_0^P E(P,\tau)d_b(\cdot,\tau)\widehat E(\tau,0)\,d\tau.
\end{equation*}
For $0<h<P/2$, define the bounded linear operator
\[
 \mathcal K_hg=\int_h^{P-h}
 E(P,\tau)\{d_b(\cdot,\tau)\widehat E(\tau,0)g\}\,d\tau\ \ \text{for  $g\in\mathcal X$.}
\]
For $\|g\|_\infty\leq1$, the support property of $d_b$, \eqref{eq4.16} and
\eqref{eq4.17} give
\[\|E(P,\tau)\{d_b(\cdot,\tau)\widehat E(\tau,0)g\}\|_\infty
\leq e^{L_f(P-\tau)}\|d_b\|_\infty e^{-2\omega {\tau}}
\leq e^{L_fP}\|d_b\|_\infty \]
and
\[
\supp\{d_b(\cdot,\tau)\widehat E(\tau,0)g\}
\subset
\{(s,y)\in\mathcal C:|s|\leq Z_0\}\ \ \text{and}\ \ 
 \|d_b(\cdot,\tau)\widehat E(\tau,0)g\|_\infty\leq\|d_b\|_\infty.
\]
Since
\[
 \left\|\nabla_X\int_{\mathcal C}H_{c_1}((t-\tau);X,Y)g(Y)\,dY\right\|_\infty
 \leq C{(t-\tau)}^{-1/2}\|g\|_\infty,
\]
then for $P-\tau\in[h,P-h]$, there holds
\[
 \begin{aligned}
 \|\nabla_XE(P,\tau)\{d_b(\cdot,\tau)\widehat E(\tau,0)g\}\|_\infty
 &\leq C\left({(P-\tau)}^{-1/2}
 +L_f\int_0^{(P-\tau)}((P-\tau)-a)^{-1/2}e^{L_fa}\,da\right)\|d_b\|_\infty\\
 &\leq C_h\|d_b\|_\infty,
 \end{aligned}
\]
where \(C_h=C\left(h^{-1/2}+2L_fe^{L_fP}\sqrt P\right).\)
For $|\xi|\leq Z_0$ and $P-\tau\leq P$, one has that
$|s+c_1(P-\tau)-\xi|\geq$ $(|s|-Z_0-c_1P)_+$.
Thus the kernel bound gives
\[
 \begin{aligned}
 |E(P,\tau)\{d_b(\cdot,\tau)\widehat E(\tau,0)g\}(X)|
 &\leq C_h\int_{-Z_0}^{Z_0}
 e^{-(s+c_1(P-\tau)-\xi)^2/(4(P-\tau))}\,d\xi\\
 &\leq C_h' e^{-(|s|-Z_0-c_1P)_+^2/(4P)}
 \end{aligned}
\]
with \(C_h=\frac{e^{L_fP}\|d_b\|_\infty}{\sqrt{4\pi h}}\) and \(C_h'=2Z_0C_h\).
Integrating these estimates in $\tau$ gives
\[
 \begin{gathered}
 \sup_{\|g\|_\infty\leq1}
 (\|\mathcal K_hg\|_\infty+\|\nabla_X\mathcal K_hg\|_\infty)<\infty,\\
 \sup_{\|g\|_\infty\leq1}\sup_{|s|\geq H}|\mathcal K_hg(s,y)|
 \leq PC_h'e^{-(H-Z_0-c_1P)_+^2/(4P)}\rightarrow 0\ \ \text{as \(H\to\infty\).}
 \end{gathered}
\]
By the Arzel\`a--Ascoli theorem
and a diagonal process, the images of
any sequence in the unit ball  have a subsequence (called \((g_j)_{j\in\mathbb N}\)) converging locally
uniformly, and it holds
\[
 \sup_{j\in\mathbb N}\sup_{|s|\geq H}|\mathcal K_hg_j|<\varepsilon/4,
 \ \ 
 \sup_{|s|\leq H}|\mathcal K_hg_j-\mathcal K_hg_k|<\varepsilon/2,
 \ \  \|\mathcal K_hg_j-\mathcal K_hg_k\|_\infty<\varepsilon.
\]
Thus $\mathcal K_h$ is compact.
Notice that 
\[
 \begin{aligned}
 &\sup_{\|g\|_\infty\leq1}
 \left\|\left(\int_0^h+\int_{P-h}^P\right)
 E(P,\tau)d_b(\cdot,\tau)\widehat E(\tau,0)g\,d\tau\right\|_\infty\\
\leq&\left(\int_0^h+\int_{P-h}^P\right)
 e^{L_f(P-\tau)}\|d_b\|_\infty e^{-2\omega {\tau}}\,d\tau
 \leq2h e^{L_fP}\|d_b\|_\infty.
 \end{aligned}
\]
Therefore
\[
 \|\mathcal M-\widehat E(P,0)-\mathcal K_h\|
 \leq2h e^{L_fP}\|d_b\|_\infty\rightarrow 0 \ \ \text{as \(h\rightarrow0\).}
\]
Hence $\mathcal M-\widehat E(P,0)$ is compact. For a bounded operator $T$,
let $\sigma(T)$ be its spectrum, $r(T)=\sup_{\lambda\in\sigma(T)}|\lambda|$
be its spectral radius, and $r_{\ess}(T)$ be the spectral radius
of its class $[T]$ in the quotient of bounded operators by compact
operators. In this quotient,
\[
 [\mathcal M]=[\widehat E(P,0)], \ \ 
 r([T])\leq\|[T]\|:=\inf_{K\in\mathcal K(\mathcal X)}\|T-K\|\leq\|T\|.
\]
Consequently
\begin{equation*}
 r_{\ess}(\mathcal M)=r_{\ess}(\widehat E(P,0))
 \leq\|\widehat E(P,0)\|\leq e^{-2\omega P}<1.
\end{equation*}
The positive function $p^0(\cdot,0)$ belongs to $\mathcal X$ and
\[
 \mathcal M p^0(\cdot,0)=E(P,0)p^0(\cdot,0)=p^0(\cdot,P)=p^0(\cdot,0),
 \ \  1\in\sigma(\mathcal M), \ \  r(\mathcal M)\geq1.
\]
Let $\mathcal X'$ be the continuous dual,
$\mathcal X'_+=\{\ell:\ell(g)\geq0\text{ for }g\geq0\}$. Define $\mathcal M'$ by
$\langle\mathcal M'\ell,g\rangle=\langle\ell,\mathcal Mg\rangle$.
Since
$g=\max\{g,0\}-\max\{-g,0\}$, Then the positive cone \(\mathcal X_+=\{g\in\mathcal X:g\geq0\}\) generates $\mathcal X$ .
Moreover,
\[
 \mathcal M\geq0, \ \  \|\mathcal M\|\leq e^{L_fP}, \ \ 
 r_{\ess}(\mathcal M)<1\leq r(\mathcal M).
\]
The cone \(\mathcal X_+\) is closed and generating, hence total.  It then follows from the  generalized
Krein--Rutman theorem \cite[Corollary 2.2, p.~324]{nussbaum-kr}
(see also \cite[Section 5.2]{gm-recovery})  that
\[
 0\ne\ell\in\mathcal X'_+, \ \ 
 \mathcal M'\ell=r(\mathcal M)\ell.
\]
Since $\mathcal C$ is locally compact, second countable and Hausdorff,
the Riesz representation theorem \cite[Theorem 1.3]{baggett-functional}
represents $\ell$ by a finite
nonzero positive Radon measure, denoted by the same letter, with
$\ell(\mathcal C):=\int_{\mathcal C}g(X)\,d\ell (X)=\|\ell\|$ for any \(g\in C_0(\mathcal C)\). Since
$[-m,m]\times\mathbb T'\uparrow\mathcal C$, there is some such compact set $A$
satisfying $\ell(A)>0$. It then follows from the continuity and positivity of $p^0$ that
\[
 0<\min_Ap^0(\cdot,0)\,\ell(A)
 \leq\langle\ell,p^0(\cdot,0)\rangle \leq \|p^0(\cdot,0)\|_\infty\,\|\ell\|<\infty.
\]
By the definition of the dual operator,
\[
 r(\mathcal M)\langle\ell,p^0(\cdot,0)\rangle
 =\langle\mathcal M'\ell,p^0(\cdot,0)\rangle
 =\langle\ell,\mathcal M p^0(\cdot,0)\rangle
 =\langle\ell,p^0(\cdot,0)\rangle.
\]
Thus
\[
 (r(\mathcal M)-1)\langle\ell,p^0(\cdot,0)\rangle=0.
\]
Hence $r(\mathcal M)=1$ and $\mathcal M'\ell=\ell$.

Rescale $\ell$ by the positive number
$\langle\ell,p^0(\cdot,0)\rangle$ so that
\[
 \langle\ell,p^0(\cdot,0)\rangle=1,
 \ \  \mathcal M'\ell=\ell.
\]
Define $\mu_t=E(P,t)'\ell$ for $0\leq t\leq P$ in the sense that \(\langle\mu_t,g\rangle
:=\langle\ell,E(P,t)g\rangle\)
for \(g\in\mathcal X\) and \(0\leq t\leq P\).
then 
\[
 \mu_0=E(P,0)'\ell=\mathcal M'\ell=\ell,
 \ \  \mu_P=E(P,P)'\ell=I'\ell=\ell.
\]
Thus $\mu_{t+kP}:=\mu_t$, for $0\leq t<P$ and $k\in\mathbb Z$,
defines a $P$-periodic family of finite positive measures,
continuous in the weak-* topology of $C_0(\mathcal C)'$.
Here a prime denotes the continuous dual operator.
The evolution law and time periodicity give
\[
 E(2P,t)=E(2P,P)E(P,t), \ \  E(2P,P)=E(P,0)=\mathcal M,
\]
hence
\[
 E(2P,t)'\ell=E(P,t)'E(2P,P)'\ell=E(P,t)'\ell=\mu_t.
\]
For $g\in C_c(\mathcal C)$, $g\geq0$, Tonelli's theorem gives
\[
 \langle\mu_t,g\rangle
 =\int_{\mathcal C}\int_{\mathcal C}
 \Gamma(2P,t;X,Y)g(Y)\,dY\,d\ell(X)
 =\int_{\mathcal C}g(Y)
 \left[\int_{\mathcal C}\Gamma(2P,t;X,Y)\,d\ell(X)\right]dY.
\]
Hence $\mu_t$ has the density
\begin{equation*}
 \psi^0(Y,t)=\int_{\mathcal C}\Gamma(2P,t;X,Y)\,d\ell(X)>0.
\end{equation*}
Thus $\mu_0=\mu_P$ implies
$\psi^0(Y,0)=\psi^0(Y,P)$ for every $Y$.
For $\tau=2P-t\in[P,2P]$, one has
\[
 \sup_{X,Y \in\mathcal C, \,0\leq t\leq P}\Gamma(2P,t;X,Y)
 \leq\frac{e^{2PL_f}}{\sqrt{4\pi P}}
 \sup_{\substack{P\leq\tau\leq2P\\y,\upsilon\in\mathbb T'}}
 H_{\mathbb T'}(\tau;y,\upsilon)=:C_P<\infty.
\]
It then follows from Tonelli's theorem and $\int H_{c_1}\,dY=1$ that 
\[
 \int_{\mathcal C}\psi^0(Y,t)\,dY
 =\int_{\mathcal C}\left[\int_{\mathcal C}\Gamma(2P,t;X,Y)\,dY\right]d\ell(X)
 \leq e^{2PL_f}\ell(\mathcal C).
\]
Consequently
\begin{equation}\label{eq4.18}
 \sup_{0\leq t\leq P}\|\psi^0(\cdot,t)\|_\infty\leq C_P\|\ell\|,
 \ \ 
 \sup_{0\leq t\leq P}\int_{\mathcal C}\psi^0(Y,t)\,dY
 \leq e^{2PL_f}\|\ell\|.
\end{equation}
For $X=({s}_X,y_X)$, $Y=({s}_Y,y_Y)$ and $|{s}_X|\leq L$, it holds 
\( |{s}_X+c_1(2P-t)-{s}_Y|\geq(|{s}_Y|-L-2c_1P)_+\).
Splitting the density integral at $|{s}_X|=L$ therefore gives
\[
 \psi^0(Y,t)\leq
 C_P\|\ell\|e^{-(|{s}_Y|-L-2c_1P)_+^2/(8P)}
 +C_P\ell(\{|{s}_X|>L\}).
\]
As a result, one obtains that 
\begin{equation}\label{eq4.19}
 \limsup_{H\to\infty}\ \sup_{|{s}_Y|\geq H,y\in\mathbb T',\ 0\leq t\leq P}
 \psi^0(Y,t)\leq C_P\ell(\{|{s}_X|>L\})\longrightarrow0
 \ \ \text{as \(L\to\infty\).}
\end{equation}

For $\varphi\in C_c^\infty(\mathcal C)$ (smooth, periodic in $y$ and
compactly supported in $s$), the heat integral equation gives
\[
 \frac{E(t+h,t)\varphi-\varphi}{h}
 \longrightarrow\Delta_X\varphi+c_1\varphi_{s}+b^0(\cdot,t)\varphi
 \quad\hbox{in }\mathcal X\ \text{as  \(h\downarrow 0\)}.
\]
For $0<t<t+h<P$, the evolution law implies
\[
 \frac{E(P,t+h)\varphi-E(P,t)\varphi}{h}
 =-E(P,t+h)\frac{E(t+h,t)\varphi-\varphi}{h}.
\]
Therefore
\[
 \begin{aligned}
 \frac{d}{dt}\int_{\mathcal C}\psi^0(X,t)\varphi(X)\,dX
 &=\left\langle\ell,\frac{d}{dt}E(P,t)\varphi\right\rangle\\
 &=-\langle\mu_t,\Delta_X\varphi+c_1\varphi_{s}+b^0(\cdot,t)\varphi\rangle\\
 &=-\int_{\mathcal C}\psi^0(X,t)
 (\Delta_X\varphi+c_1\varphi_{s}+b^0(X,t)\varphi)\,dX.
 \end{aligned}
\]
Hence
\begin{equation}\label{eq4.20}
 -\psi_t^0=\Delta_X\psi^0-c_1\psi_{s}^0+b^0\psi^0.
\end{equation}
Set $\widetilde\psi(X,\tau)=\psi^0(X,-\tau)$. Then there holds
\[
 (\partial_\tau-\Delta_X+c_1\partial_{s})\widetilde\psi
 =b^0(X,-\tau)\widetilde\psi.
\]
By the interior $L^q$ estimate
\cite[Theorem 7.22]{lieberman-parabolic}, the
Sobolev embedding 
and the interior Schauder regularity theorem \cite[Theorem 4.9]{lieberman-parabolic}, one obtains that 
\(\widetilde\psi,\ \psi^0\in C^{2+\alpha/2,1+\alpha/4}_{\mathrm{loc}}
(\mathcal C\times\mathbb R)\).
Recall that $p^0(\cdot,t)\in\mathcal X$. Uniqueness for the linearized equation
and the evolution law yield that
\[
 p^0(\cdot,t)=E(t,0)p^0(\cdot,0)\ \ \text{and}\ \ 
 E(P,t)p^0(\cdot,t)=p^0(\cdot,P)=p^0(\cdot,0)\ \ \text{for $0\leq t\leq P$}.
\]
By the density representation and the definition of  the dual operator, one arrives at 
\[
 \begin{aligned}
 \int_{\mathcal C}p^0(X,t)\psi^0(X,t)\,dX
 =\langle\mu_t,p^0(\cdot,t)\rangle=\langle\ell,E(P,t)p^0(\cdot,t)\rangle=\langle\ell,p^0(\cdot,0)\rangle=1.
 \end{aligned}
\]
Periodicity extends this equality to every $t\in\mathbb R$.
For a constant $r$,  since \((s-c_1r)+c_1(t+r)=s+c_1t\), direct calculations imply that
\(V^r(s,y,t) =V^0(s-c_1r,y,t+r)\), \( p^r(s,y,t) =p^0(s-c_1r,y,t+r)\) and \(b^r(s,y,t) =b^0(s-c_1r,y,t+r)\).
It then follows from \eqref{eq4.20} that
\[
 -\psi_t^r(s,y,t)=\Delta_X\psi^r-c_1\psi_{s}^r+b^r\psi^r.
\]
Let $\xi=s-c_1r$. Then  $d\xi=ds$, gives
\[
 \begin{aligned}
 \int_{\mathcal C}p^r(X,t)\psi^r(X,t)\,dX
 &=\int_{\mathbb R\times\mathbb T'}
 p^0(\xi,y,t+r)\psi^0(\xi,y,t+r)\,d\xi\,dy=1.
 \end{aligned}
\]
This proves \eqref{eq4.12}. The proof is complete.
\end{pr}

At fixed $r$, it follows from \eqref{eq4.12} and \eqref{eq4.13} that
\begin{align*}
 \partial_t(p^r\psi^r)
 =\psi^r\Delta_Xp^r-p^r\Delta_X\psi^r+c_1\partial_{s}(p^r\psi^r)=\diver_X(c_1e_{s}p^r\psi^r+2\psi^r\nabla_Xp^r)
 -\Delta_X(p^r\psi^r).
\end{align*}
namely,
\begin{equation}\label{eq4.21}
 \partial_t(p^r\psi^r)+\Delta_X(p^r\psi^r)
 -\diver_X\bigl(c_1e_{s}p^r\psi^r+2\psi^r\nabla_Xp^r\bigr)=0
\end{equation}

\begin{lem}\label{lem-t17-weight-estimates}
There are constants $C,\eta>0$, independent of $r,X,t$, such that
\begin{equation}\label{eq4.22}
 \begin{gathered}
 p^r+\psi^r\leq C e^{-\eta|s-c_1r|},\\
 |\nabla_X\log p^r|+|\nabla_X\log\psi^r|
 +|\partial_r\log p^r|+|\partial_r\log\psi^r|\leq C
 , \ \   r\in\mathbb R.
 \end{gathered}
\end{equation}
\end{lem}

\begin{pr}
Let $\widetilde\psi(X,t)=\psi^0(X,-t)$ be given as  in Lemma \ref{lem-t17-weight}. 
Then one gets that
\[
 \begin{aligned}
 p_t^0-\Delta_Xp^0-c_1p_{s}^0-b^0(X,t)p^0&=0,\\
 \widetilde\psi_t-\Delta_X\widetilde\psi
 +c_1\widetilde\psi_{s}-b^0(X,-t)\widetilde\psi&=0.
 \end{aligned}
\]
Moreover, one infers from \eqref{eq4.2} and \eqref{eq4.18} that
\[
 0<p^0(X,t)\leq M, \ \ 
 0<\widetilde\psi(X,t)\leq C_P\|\ell\|.
\]
By \eqref{eq4.14} and \eqref{eq4.19}, one infers that 
\[
 \lim_{H\to\infty}
 \sup_{\substack{|s|\geq H,\ y\in\mathbb T',\ 0\leq t\leq P}}
 \bigl(p^0(s,y,t)+\widetilde\psi(s,y,t)\bigr)=0.
\]
By \eqref{eq4.15}, for every $t\in\mathbb R$,
\[
 b^0(X,t)\leq-4\omega\ \ \text{and}\ \ b^0(X,-t)\leq-4\omega
 \ \ \text{for \(|s|\geq Z_0\)}.
\]
Choose
\[
 0<\eta<\min\left\{1,\frac{2\omega}{1+c_1}\right\}.
\]
Then $\eta^2+c_1\eta<(1+c_1)\eta<2\omega$. Define
\[
 C_\pm= \max\bigg\{\max_{\substack{y\in\mathbb T',\ 0\leq t\leq P}}p^0(\pm Z_0,y,t),
\max_{\substack{y\in\mathbb T',\ 0\leq t\leq P}}\widetilde\psi(\pm Z_0,y,t)\bigg\}>0.
\]
Set
\[
 g_+(s)=C_+e^{-\eta(s-Z_0)}\ \ \text{for \(s\geq Z_0\) \ \ and }\ \ 
 g_-(s)=C_-e^{\eta(s+Z_0)}\ \ \text{for \(s\leq-Z_0\)}.
\]
It is obvious that 
\[
 p^0(\pm Z_0,y,t)\leq g_\pm(\pm Z_0)\ \ \text{and}\ \ 
 \widetilde\psi(\pm Z_0,y,t)\leq g_\pm(\pm Z_0)\ \ \text{for all $(y,t)\in \mathbb{T}'\times [0,P]$}.
\]
Since $g_\pm$ depend only on $s$, it then can be inferred that 
\(\partial_tg_\pm=0\),  \(\partial_{s}g_\pm=\mp\eta g_\pm\) and \(\Delta_Xg_\pm=\partial_{ss}g_\pm=\eta^2g_\pm. \)
Thus
\[
 \begin{aligned}
 (\partial_t-\Delta_X-c_1\partial_{s}-b^0(X,t))g_\pm
 &=(-\eta^2\pm c_1\eta-b^0(X,t))g_\pm\\
 &\geq(4\omega-\eta^2-c_1\eta)g_\pm>2\omega g_\pm>0,\ \ s\geq Z_0\\
 (\partial_t-\Delta_X+c_1\partial_{s}-b^0(X,-t))g_\pm
 &=(-\eta^2\mp c_1\eta-b^0(X,-t))g_\pm\\
 &\geq(4\omega-\eta^2-c_1\eta)g_\pm>2\omega g_\pm>0,\ \ s\leq -Z_0.
 \end{aligned}
\]

We prove $p^0\leq g_+$ on $s\geq Z_0$, $y\in\mathbb T'$ and $t\in[0,P]$. 
If this is not true, there would be
${s}_1>Z_0$, $y_1\in\mathbb T'$, and $t_1\in[0,P]$ such that
\[
 p^0({s}_1,y_1,t_1)-g_+({s}_1)>0.
\]
Since  $p^0$ and $g_+$ vanish uniformly as $s\to+\infty$. Hence there is some
$H>{s}_1$ such that 
\[
 \sup_{\substack{s\geq H,\ y\in\mathbb T',\ 0\leq t\leq P}}
 |p^0(s,y,t)-g_+(s)|
 <\frac12\bigl(p^0({s}_1,y_1,t_1)-g_+({s}_1)\bigr).
\]
On the compact cylinder $[Z_0,H]\times\mathbb T'\times[0,P]$, the
continuous function $p^0-g_+$ attains a maximum at some
$({s}_*,y_*,t_*)$, with \(Z_0<{s}_*<H\) and
\[
 \begin{aligned}
 p^0({s}_*,y_*,t_*)-g_+({s}_*)
 \geq p^0({s}_1,y_1,t_1)-g_+({s}_1)>0
 \end{aligned}
\]
and
\[
 \partial_t(p^0-g_+)\geq 0, \ \ 
 \nabla_X(p^0-g_+)=0, \ \ 
 \Delta_X(p^0-g_+)\leq0.
\]
Consequently, at $({s}_*,y_*,t_*)$,
\[
 \begin{aligned}
 (\partial_t-\Delta_X-c_1\partial_{s}-b^0)(p^0-g_+)\geq-b^0(p^0-g_+)
 \geq4\omega(p^0-g_+)>0.
 \end{aligned}
\]
But 
\[
 \begin{aligned}
 (\partial_t-\Delta_X-c_1\partial_{s}-b^0)(p^0-g_+)
 =-(\partial_t-\Delta_X-c_1\partial_{s}-b^0)g_+<0.
 \end{aligned}
\]
A contradiction. Thus $p^0\leq g_+$ on $s\geq Z_0$, $y\in\mathbb T'$ and $t\in[0,P]$. 
Similarly, one has that $p^0\leq g_-$ on \(s\leq -Z_0\), $y\in\mathbb T'$ and $t\in[0,P]$. 
Also there holds \(\widetilde\psi(s,y,t)\leq g_+(s)\) 
on $s\geq Z_0$, $y\in\mathbb T'$ and $t\in[0,P]$, and 
\(\widetilde\psi(s,y,t)\leq g_-(s)\) 
on \(s\leq -Z_0\), $y\in\mathbb T'$ and $t\in[0,P]$.  
Choose a constant $C$ so that
\[
 C\geq e^{\eta Z_0}\max\left\{
 2C_+,\ 2C_-,\
 \max_{\substack{|s|\leq Z_0,\ y\in\mathbb T',\ 0\leq t\leq P}}
 (p^0(s,y,t)+\psi^0(s,y,t))\right\}.
\]
It can be inferred that 
\[
 p^0+\psi^0\leq2C_+e^{\eta Z_0}e^{-\eta|s|}
   \leq Ce^{-\eta|s|}\ \ \text{on \(s\geq Z_0\), $y\in\mathbb T'$ and $t\in[0,P]$}
\]
and
\[
 p^0+\psi^0\leq2C_-e^{\eta Z_0}e^{-\eta|s|}
   \leq Ce^{-\eta|s|}\ \ \text{on \(s\leq -Z_0\), $y\in\mathbb T'$ and $t\in[0,P]$. }
\]
For $|s|\leq Z_0$, the choice of $C$ yields that 
\[
 p^0(X,t)+\psi^0(X,t)
 \leq Ce^{-\eta Z_0}\leq Ce^{-\eta|s|}.
\]
Combining these three inequalities, one infers that 
\begin{equation}\label{eq4.23}
 p^0(X,t)+\psi^0(X,t)\leq Ce^{-\eta|s|}
\ \ \text{for \(X\in\mathcal C,\ t\in\mathbb R\)}.
\end{equation}

By the interior Schauder estimate, the parabolic Harnack
inequality and $P$-periodicity in time, there exists a
constant $C>0$, independent of $X,t$, such that
\[
 \begin{aligned}
 |\nabla_Xp^0(X,t)|+|p_t^0(X,t)|
 &\leq Cp^0(X,t),\\
 |\nabla_X\widetilde\psi(X,t)|+|\widetilde\psi_t(X,t)|
 &\leq C\widetilde\psi(X,t),
 \end{aligned}
 \ \  (X,t)\in\mathcal C\times\mathbb R.
\]
Since $\psi^0(X,t)=\widetilde\psi(X,-t)$, one has
\[
 \begin{aligned}
 |\nabla_X\psi^0(X,t)|+|\psi_t^0(X,t)|
 =|\nabla_X\widetilde\psi(X,-t)|
   +|\widetilde\psi_t(X,-t)|\leq C\widetilde\psi(X,-t)
 =C\psi^0(X,t).
 \end{aligned}
\]
Dividing the estimates for $p^0$ and $\psi^0$ by these
strictly positive functions and enlarging $C$ if necessary,
we obtain
\begin{equation}\label{eq4.24}
 |\nabla_X\log p^0|
 +\frac{|p_t^0|}{p^0}
 +|\nabla_X\log\psi^0|
 +\frac{|\psi_t^0|}{\psi^0}
 \leq C
 \ \ \text{in }\mathcal C\times\mathbb R.
\end{equation}
By the definitions of $p^r$ and $\psi^r$, we have
\[
 \begin{aligned}
 p^r(s,y,t)&=p^0(s-c_1r,y,t+r),\\
 \psi^r(s,y,t)&=\psi^0(s-c_1r,y,t+r),
 \end{aligned}
 \ \  (s,y,t)\in\mathcal C\times\mathbb R,\quad r\in\mathbb R.
\]
It follows from \eqref{eq4.23} that   
\[
 p^r(X,t)+\psi^r(X,t)\leq Ce^{-\eta|s-c_1r|}.
\]
After some calculations, one has that 
\[
 \begin{gathered}
 \nabla_X\log p^r(s,y,t)
    =(\nabla_X\log p^0),\ \ 
 \nabla_X\log\psi^r(s,y,t)
    =(\nabla_X\log\psi^0),\\
 \partial_r\log p^r(s,y,t)
 =\frac{-c_1p_{s}^0+p_t^0}{p^0},\ \
 \partial_r\log\psi^r(s,y,t)
 =\frac{-c_1\psi_{s}^0+\psi_t^0}{\psi^0}, 
 \end{gathered}
\] where \(p^0\) and \(\psi^0\) are evaluated at \((s-c_1r,y,t+r)\).
It then can be derived from \eqref{eq4.24} that
\[
 \begin{aligned}
 |\nabla_X\log p^r(X,t)|
 +|\nabla_X\log\psi^r(X,t)|
 &\leq C,\\
 |\partial_r\log p^r(X,t)|
 +|\partial_r\log\psi^r(X,t)|
 &\leq (1+c_1)C,
 \end{aligned}
 \ \  (X,t)\in\mathcal C\times\mathbb R,\quad r\in\mathbb R,
\]
where $C$ is independent of $X,t,r$.
This proves \eqref{eq4.22}.
The proof is complete.
\end{pr}
\vspace{0.3cm}

For $x\in\mathbb R^N$, since $W_t>0$ and $W(x,-\infty)=0$, $W(x,+\infty)=1$,
then for every $v\in(0,1)$, there exists a unique
$\tau(x,v)\in\mathbb R$ such that \( W(x,\tau(x,v))=v\).
Define $J$ by
\begin{equation*}
 J(x,v)=W_t(x,\tau(x,v)),\ \ 0<v<1,\ \ 
 J(x,v)=0,\ \ v\notin(0,1).
\end{equation*}

\begin{lem}\label{lem-t17-J}
The function $J$ is periodic in $x$, globally Lipschitz in $(x,v)$,
and satisfies
\begin{equation}\label{eq4.25}
 0\leq J\leq M, \ \ 
 J(x,v)\geq {\kappa}\ \ \text{for   \(\eta_0\leq v\leq1- \eta_0\)}, \ \ 
 J(x,W^r(x,t))=\partial_rW^r(x,t),
\end{equation}
where \(M\), \(\eta_0\) and \(\kappa\)  are given in \eqref{eq4.2}.
\end{lem}
\begin{pr}
Since $W_t-\Delta_xW=f(x,W)$, $0<W<1$, and $f$ is uniformly
bounded on $\mathbb R^N\times[0,1]$, the interior parabolic
gradient estimate  implies that 
\[
 |\nabla_xW(x,\tau(x,v))|\leq C,
 \ \  x\in\mathbb R^N,\quad 0<v<1.
\]
At $(x,\tau(x,v))$, there hold
\[
 \tau_v=W_t^{-1}, \ \  \tau_{x_i}=-W_{x_i}/W_t,
 \ \  J_v=W_{tt}/W_t,
 \ \  J_{x_i}=W_{tx_i}-(W_{tt}/W_t)W_{x_i}.
\]
In moving coordinates $X=(x_1-c_1t,y)$, one infers that
\[
 W_t=p^0(X,t), \ \  W_{tt}=p_t^0(X,t)-c_1p_{s}^0(X,t),
 \ \  \nabla_xW_t=\nabla_Xp^0(X,t).
\]
By Lemma \ref{lem-t17-weight-estimates}, one has that 
\begin{equation}\label{eq4.26}
 |J_v|+|\nabla_xJ|\leq C,\ \ 0<v<1.
\end{equation}
For every $H>0$, it follows from the periodicity and strict positivity that
\[
 \min_{x\in\mathbb{R}^N} U_1(x,H)>0\ \ \text{and} \ \ \min_{x\in\mathbb{R}^N}(1-U_1(x,-H))>0.
\]
Hence $v\downarrow0$ forces $x_1-c_1\tau(x,v)\to+\infty$, and
$v\uparrow1$ forces $x_1-c_1\tau(x,v)\to-\infty$, uniformly in $x$.
Recall that \(W_t(x,t)=p^0(x_1-c_1t,y,t)\) for \(x=(x_1,y)\).
It follows from \eqref{eq4.23} that
\[
 \lim_{v\downarrow0}\sup_{x\in\mathbb{R}^N}J(x,v)
 =\lim_{v\uparrow1}\sup_{x\in\mathbb{R}^N}J(x,v)=0.
\]
The mean value theorem  and \eqref{eq4.26}  yield that there is a constant  $C_J$ such that 
\begin{equation*}
 |J(x,v)-J(\widetilde x,\widetilde v)|
 \leq C_J(|x-\widetilde x|+|v-\widetilde v|)
 ,\ \ x,\widetilde x\in\mathbb R^N,\ v,\widetilde v\in\mathbb R.
\end{equation*}
For a period vector $kL$, there hold
\[
 W(x+kL,t)=W(x,t-k_1L_1/c_1),\quad
 \tau(x+kL,v)=\tau(x,v)+k_1L_1/c_1,
 \quad J(x+kL,v)=J(x,v).
\]
The other assertions in \eqref{eq4.25} follow from
\eqref{eq4.2} and the inverse identity. The proof is complete.
\end{pr}
\vspace{0.3cm}

Let $D$ be fixed as in Lemma \ref{lem-t17-trapping}, $T\in\mathbb R$,
and $\sigma\in\{-1,1\}$. Assume
\begin{equation}\label{eq4.27}
 \begin{gathered}
 a\in C^2([T,\infty)),\quad a(T)=0,\quad a'\geq0,\quad
 a(t)\to a_\infty<\infty,\quad a'(t)\to0,\\
 \|a'\|_\infty+\|a''\|_\infty<\infty,
 \ \  h_0\in C_b^2(\mathbb R^N),\quad0\leq h_0\leq D.
 \end{gathered}
\end{equation}
Here $C_b^2$ means that the function and its derivatives of order at
most two are bounded.
For each $\sigma\in\{-1,1\}$, consider the initial-value problem
\begin{equation}\label{eq4.28}
 \begin{cases}
 (Z_{\sigma})_t=\Delta Z_{\sigma}+f(x,Z_{\sigma})+\sigma a'(t)J(x,Z_{\sigma}),\ \ 
x\in\mathbb R^N,\ t>T,\\
 Z_{\sigma}(x,T)=W^{\sigma h_0(x_1-c_1T,y)}(x,T),\ \ x\in\mathbb R^N.
 \end{cases}
\end{equation}
We write $Z_+=Z_1$ and $Z_-=Z_{-1}$ for simplicity.

\begin{lem}\label{lem-t17-auxiliary}
Under \eqref{eq4.27}, problem \eqref{eq4.28} has a unique classical 
bounded solution  for $t>T$, and
\begin{equation}\label{eq4.29}
 W^{a(t)}\leq Z_+\leq W^{a(t)+D}, \ \ 
 W^{-a(t)-D}\leq Z_-\leq W^{-a(t)}.
\end{equation}
\end{lem}
\begin{pr}
By Lemma \ref{lem-t17-J}, the reaction in \eqref{eq4.28}
is continuous in time and globally Lipschitz in the dependent
variable, uniformly for $t\geq T$.
Thus \cite[Chapter 6, Theorem 1.2]{pazy}, applied on
$BUC(\mathbb R^N)$, gives a unique global mild solution $Z_\sigma$.
Since $f(x,0)=f(x,1)=J(x,0)=J(x,1)=0$, the comparison principle
yields $0\leq Z_\sigma\leq1$.
Interior parabolic regularity and the Schauder estimate then
imply that $Z_\sigma$ is a classical solution for $t>T$.
For every constant $d\in[0,D]$, there holds
\[
 \Nop[W^{\sigma(a(t)+d)}]
 =\sigma a'(t)\partial_rW^r|_{r=\sigma(a+d)}
 =\sigma a'(t)J(x,W^{\sigma(a+d)}).
\]
At $t=T$, one has that
\[
 W\leq Z_+(T)\leq W^D,
 \ \  W^{-D}\leq Z_-(T)\leq W.
\]
It then follows from the comparison principle that \eqref{eq4.29} holds.
The proof is complete.
\end{pr}

\begin{lem}\label{lem-t17-phase}
Under \eqref{eq4.27}, for each sign there is a unique phase
$h=h_{\sigma}$ such that
\begin{equation}\label{eq4.30}
 Z_{\sigma}((s+c_1t,y),t)=V^{\sigma(a(t)+h(X,t))}(X,t),
 \ \ 0\leq h\leq D, \ \  h(X,T)=h_0(X),
\end{equation}
It is classical locally for $t>T$ and satisfies
\begin{equation}\label{eq4.31}
 h_t=\Delta_Xh+
 \left.(c_1e_{s}+2\nabla_X\log p^r)\right|_{r=\sigma(a+h)}\cdot\nabla_Xh
 +\sigma\left.\partial_r\log p^r\right|_{r=\sigma(a+h)}|\nabla_Xh|^2.
\end{equation}
\end{lem}
\begin{pr}
Strict monotonicity in $r$ and \eqref{eq4.29} define
the unique $h$ in \eqref{eq4.30}. For \(H>0\), one has that 
\[
 \min_{\substack{|s|\leq H,\ |r|\leq D+a_\infty\\
y\in\mathbb T',\ t\in[0,P]}}p^r(X,t)>0.
\]
Since
\[
 \partial_h V^{\sigma(a(t)+h)}(X,t)
 =\sigma p^{\sigma(a(t)+h)}(X,t)\ne0,
\]
the implicit function theorem, together with the classical
regularity of $Z_\sigma$ and the regularity of $V^r$ and $a$,
yields
\[
 h\in C_{\mathrm{loc}}^{2,1}
 \bigl(\mathbb R^N\times(T,\infty)\bigr).
\]
Differentiating
\eqref{eq4.30} with all derivatives of $V^r,p^r$ at fixed $r$,
gives
\[
 \begin{aligned}
 \partial_t[Z_{\sigma}((s+c_1t,y),t)]
 &=\left.V_t^r\right|_{r=\sigma(a+h)}
+\sigma p^{\sigma(a+h)}(a'+h_t),\\
 \partial_s[Z_\sigma((s+c_1t,y),t)]
 &=\left.V_s^r\right|_{r=\sigma(a+h)}
   +\sigma p^{\sigma(a+h)}h_s,\\
 \Delta_X[Z_{\sigma}((s+c_1t,y),t)]
 &=\left.\Delta_XV^r\right|_{r=\sigma(a+h)}
+2\sigma\left.\nabla_Xp^r\right|_{r=\sigma(a+h)}\cdot\nabla_Xh\\
 &\quad+\left.\partial_rp^r\right|_{r=\sigma(a+h)}|\nabla_Xh|^2
+\sigma p^{\sigma(a+h)}\Delta_Xh.
 \end{aligned}
\]
The equation \eqref{eq4.28} in moving coordinates is
\[
 \begin{aligned}
 \partial_t[Z_{\sigma}((s+c_1t,y),t)]
 &=(\Delta_X+c_1\partial_{s})[Z_{\sigma}((s+c_1t,y),t)]\\
 &\quad+f((s+c_1t,y),Z_{\sigma}((s+c_1t,y),t))+\sigma a'p^{\sigma(a+h)}.
 \end{aligned}
\]
Since \( h_s=e_s\cdot\nabla_Xh\), where \(e_s=(1,0,\ldots,0)\), it follows from \eqref{eq4.13} that
\[
 \begin{aligned}
 {\sigma} p^{\sigma(a+h)} h_t
 &=\sigma p^{\sigma(a+h)}\Delta_Xh
  +\sigma\left.(c_1p^re_{s}+2\nabla_Xp^r)\right|_{r=\sigma(a+h)}\cdot\nabla_Xh
+\left.\partial_rp^r\right|_{r=\sigma(a+h)}|\nabla_Xh|^2.
 \end{aligned}
\]
Since $1/\sigma=\sigma$, one has 
\[
 h_t=\Delta_Xh+
\left.(c_1e_{s}+2\nabla_X\log p^r)\right|_{r=\sigma(a+h)}
\cdot\nabla_Xh
+\sigma\left.\partial_r\log p^r\right|_{r=\sigma(a+h)}
|\nabla_Xh|^2.
\]
This completes the proof.
\end{pr}

Use the setting of Lemma \ref{lem-t17-phase} and extend $p^r,\psi^r$
periodically in $y$ to $\mathbb R^N$. Fix $\sigma\in\{-1,1\}$ and let
$h=h_{\sigma}$ be the phase in \eqref{eq4.30}. Define
\begin{equation*}
 \begin{gathered}
 F(X,t)=\int_0^{h(X,t)}m(X,t,q)\,dq,\ \ 0\leq q\leq D,
 \end{gathered}
\end{equation*}
where \(m(X,t,q)=e^{K_0q}p^{\sigma(a(t)+q)}(X,t)\psi^{\sigma(a(t)+q)}(X,t)\)
and \( K_0=1+\sup_{\substack{r,t\in\mathbb R,\ X\in\mathcal C}}
 |\partial_r\log\psi^r(X,t)|\)
is finite by \eqref{eq4.22}.

\begin{lem}\label{lem-t17-entropy}
Under the hypotheses of Lemma \ref{lem-t17-phase}, assume additionally
that $h_0(s,y)=0$ for $y$ outside a fixed compact subset of
$\mathbb R^{N-1}$.
Then, there is a constant $C$ independent of $X,t$ such that
\begin{equation}\label{eq4.32}
 \sup_{t\geq T}\int_{\mathbb R^N}F(X,t)\,dX
 +\int_T^\infty\!\int_{\mathbb R^N}
 m(X,t,h(X,t))|\nabla_Xh|^2\,dX\,dt
 \leq2e^{Ca_\infty}\int_{\mathbb R^N}F(X,T)\,dX<\infty.
\end{equation}
\end{lem}
\begin{pr}
It can be inferred from the identity \eqref{eq4.21} that
\[
 \begin{aligned}
 m_t&=-\Delta_Xm
 +\diver_X\big(m(c_1e_{s}+2\nabla_X\log p^{\sigma(a+q)})\big)
 +\sigma a'm\,\partial_r\log(p^r\psi^r)|_{r=\sigma(a+q)},\\
 \partial_qm&=m\big(K_0+\sigma\partial_r\log(p^r\psi^r)|_{r=\sigma(a+q)}\big).
 \end{aligned}
\]
Define 
\[
 Q(X,t)=\int_0^h
 m\big(c_1e_{s}-2\nabla_X\log\psi^{\sigma(a+q)}\big)\,dq
 =\int_0^h\big(m(c_1e_{s}+2\nabla_X\log p^{\sigma(a+q)})-2\nabla_Xm\big)\,dq.
\]
Then
\[
 \begin{aligned}
 F_t&=\int_0^h m_t\,dq+m(X,t,h(X,t))h_t,\\
 \Delta_XF&=\int_0^h\Delta_Xm\,dq
 +2(\nabla_Xm)(X,t,h(X,t))\cdot\nabla_Xh\\
 &\quad+(\partial_qm)(X,t,h(X,t))|\nabla_Xh|^2
 +m(X,t,h(X,t))\Delta_Xh,\\
 \diver_XQ&=\int_0^h\diver_X\big(m(c_1e_{s}+2\nabla_X\log p^{\sigma(a+q)})
  -2\nabla_Xm\big)\,dq\\
 &\quad+\big(m(X,t,h(X,t))(c_1e_{s}+2\nabla_X\log p^{\sigma(a+h)})
     -2(\nabla_Xm)(X,t,h(X,t))\big)\cdot\nabla_Xh.
 \end{aligned}
\]
It then follows from  \eqref{eq4.31} that
\begin{equation}\label{eq4.33}
 \begin{aligned}
 F_t={}&\Delta_XF+\diver_XQ
 -m(X,t,h(X,t))\big(K_0+\sigma\partial_r\log\psi^r|_{r=\sigma(a+h)}\big)|\nabla_Xh|^2\\
 &+\sigma a'\int_0^h m\,\partial_r\log(p^r\psi^r)|_{r=\sigma(a+q)}\,dq.
 \end{aligned}
\end{equation}
For each $H>0$, define
\[
 c_H:=
 \min_{\substack{|s|\leq H,\ |r|\leq D+a_\infty\\
                  y\in\mathbb T',\ t\in[0,P]}}
 p^r(s,y,t)\psi^r(s,y,t)>0.
\]
By the choice of $K_0$ and Lemma \ref{lem-t17-weight-estimates}, one obtains that 
\begin{equation}\label{eq4.34}
 \begin{gathered}
 K_0+\sigma\partial_r\log\psi^r\geq1, \ \  |Q|\leq CF, \ \ 
 \left|\int_0^h m\,\partial_r\log(p^r\psi^r)|_{r=\sigma(a+q)}\,dq\right|\leq CF,\\
 0\leq F\leq Ce^{-\eta|s|}, \ \ 
 F\geq c_Hh, \ \  m(X,t,h(X,t))\geq c_H>0,\ \ |s|\leq H.
 \end{gathered}
\end{equation}
The constants may depend on $D+a_\infty$.
Recall that \(h(X,T)=h_0(X)\) from \eqref{eq4.30}. By the assumption of the lemma, there is a compact set 
\(C_0\subset\mathbb R^{N-1}\) such that \(h_0(s,y)=0\) for \(y\notin C_0\). Hence \(F(s,y,T)=0\) for \(y\notin C_0\).
Hence
\[
 \int_{\mathbb R^N}F(X,T)\,dX
 \leq C\int_{\mathbb R}e^{-\eta|s|}\,ds\,|C_0|
 =\frac{2C}{\eta}|C_0|<\infty.
\]

For $n\geq1$, define \(\chi_n(y)=e^{-\sqrt{1+|y|^2}/n}\). Then there is a constant \(C\)
such that \( |\nabla\chi_n|+|\Delta\chi_n|\leq(C/n)\chi_n.\) For \(H,L>0\), 
choose smooth cutoffs $\theta_H(s),\xi_L(y)$, equal to $1$ on
$|s|\leq H$, $|y|\leq L$, and zero outside $|s|<2H$, $|y|<2L$,
with values in $[0,1]$, then there hold
\[
 |\theta_H'|\leq C/H,\quad|\theta_H''|\leq C/H^2,
 \ \  |\nabla\xi_L|\leq C/L,\quad|\Delta\xi_L|\leq C/L^2.
\]
Multiplying \eqref{eq4.33} by
$\theta_H(s)\xi_L(y)\chi_n(y)$, integrating over $[T+\nu,S]\times \mathbb{R}^N$,
and using \eqref{eq4.34}, one obtains that
\[
 \begin{aligned}
 &\int_{\mathbb R^N} \theta_H\xi_L\chi_n F(X,S)\,dX
   +\int_{T+\nu}^S\!\int_{\mathbb R^N} \theta_H\xi_L\chi_n
    m(X,t,h(X,t))|\nabla_X h|^2\,dX\,dt\\
 &\quad\leq\int_{\mathbb R^N} \theta_H\xi_L\chi_n F(X,T+\nu)\,dX\\
 & \ \ +\int_{T+\nu}^S\!\int_{\mathbb R^N}
 \{F\Delta_X(\theta_H\xi_L\chi_n)
   -Q\cdot\nabla_X(\theta_H\xi_L\chi_n)\}\,dX\,dt\\
 & \ \ +C\int_{T+\nu}^S a'(t)
 \int_{\mathbb R^N}\theta_H\xi_L\chi_n F(X,t)\,dX\,dt.
 \end{aligned}
\]
Here $S>T+\nu>T$. Note that
\[
\begin{aligned}
&F\Delta_X(\theta_H\xi_L\chi_n)
-Q\cdot\nabla_X(\theta_H\xi_L\chi_n)\\
={}&
F\theta_H''\xi_L\chi_n
-\theta_H'\xi_L\chi_n(Q\cdot e_{s})\\
&+\theta_HF\chi_n\Delta_y\xi_L
+2\theta_HF\,\nabla_y\xi_L\cdot\nabla_y\chi_n
-\theta_H\chi_n\,Q\cdot(0,\nabla_y\xi_L)\\
&+\theta_HF\xi_L\Delta_y\chi_n
-\theta_H\xi_L\,Q\cdot(0,\nabla_y\chi_n).
\end{aligned}
\]
For fixed $n,L,S$, there holds
\[\begin{aligned}
&F\theta_H''\xi_L\chi_n
-(Q\cdot e_{s})\theta_H'\xi_L\chi_n 
\leq    C(H^{-1}+H^{-2})
 \int_T^S\!\int_{|s|\geq H}\!\int_{|y|\leq2L}
 e^{-\eta|s|}\chi_n(y)\,dy\,ds\,dt\rightarrow0
\end{aligned}
\] as \(H\to\infty\).
After $H\to\infty$,  it holds
\[\begin{aligned}
&F\chi_n\Delta_y\xi_L
+2F\nabla_y\xi_L\cdot\nabla_y\chi_n
-\chi_nQ\cdot(0,\nabla_y\xi_L)\\
\leq& 
 C(L^{-1}+L^{-2}+L^{-1}n^{-1})
 \int_T^S\!\int_{\mathbb R^N}e^{-\eta|s|}\chi_n(y)\,dX\,dt
 \rightarrow0\ \ \text{as \(L\to\infty\)}.
\end{aligned}
\]
It then follows from dominated convergence and Fatou's
lemma that
\[
 \begin{aligned}
 &\int_{\mathbb R^N}\chi_n(y)F(X,S)\,dX
 +\int_T^S\!\int_{\mathbb R^N}
    \chi_n(y)m(X,t,h(X,t))|\nabla h|^2\,dX\,dt\\
 \leq&\int_{\mathbb R^N}\chi_n(y)F(X,T)\,dX+\int_T^S(C/n+Ca'(t))
   \int_{\mathbb R^N}\chi_n(y)F(X,t)\,dX\,dt.
 \end{aligned}
\]
By Gronwall's inequality, one arrives at
\[
 \int_{\mathbb R^N}\chi_n(y)F(X,t)\,dX
 \leq e^{C(t-T)/n+Ca(t)}
\int_{\mathbb R^N}\chi_n(y)F(X,T)\,dX.
\]
Hence, it holds
\[
 \begin{aligned}
 &\int_{\mathbb R^N}\chi_n(y)F(X,S)\,dX
 +\int_T^S\!\int_{\mathbb R^N}
     \chi_n(y)m(X,t,h(X,t))|\nabla h|^2\,dX\,dt\\
 \leq &e^{C(S-T)/n+Ca(S)}
 \int_{\mathbb R^N}\chi_n(y)F(X,T)\,dX\\
\leq  &e^{C(S-T)/n+Ca_\infty}
 \int_{\mathbb R^N}F(X,T)\,dX.
 \end{aligned}
\]
For each finite $S$,  let $n\to\infty$, one has 
\[
\begin{aligned}
\int_{\mathbb R^N}F(X,S)\,dX
+\int_T^S\!\int_{\mathbb R^N}
m(X,t,h)|\nabla_Xh|^2\,dX\,dt\leq e^{Ca_\infty}
\int_{\mathbb R^N}F(X,T)\,dX.
\end{aligned}
\]
Taking $S\rightarrow +\infty $, one gets 
\[
\begin{aligned}
\sup_{t\geq T}\int_{\mathbb R^N}F(X,t)\,dX
+\int_T^\infty\!\int_{\mathbb R^N}
m(X,t,h)|\nabla_Xh|^2\,dX\,dt \leq
2e^{Ca_\infty}\int_{\mathbb R^N}F(X,T)\,dX<\infty.
\end{aligned}
\]
The proof is complete.
\end{pr}

\begin{lem}\label{lem-t17-decay}
Under the hypotheses of Lemma \ref{lem-t17-entropy}, for each $H>0$, it holds
\begin{equation}\label{eq4.35}
 \sup_{|s|\leq H,\ y\in\mathbb R^{N-1}}h_{\sigma}(s,y,t)\longrightarrow0\ \ \text{as \(t\rightarrow+\infty\)}.
\end{equation}
Moreover, there holds
\begin{equation}\label{eq4.36}
 \|Z_{\sigma}(\cdot,t)-W^{\sigma a_\infty}(\cdot,t)\|_{L^\infty(\mathbb R^N)}
 \longrightarrow0,\ \ \sigma=\pm1.
\end{equation}
\end{lem}
\begin{pr}
Fix a sign and keep the phase $h$ from Lemma \ref{lem-t17-phase}.
By \eqref{eq4.32},
\begin{equation}\label{eq4.37}
 \begin{aligned}
& \int_{\mathbb R^N}F(X,t)\,dX
 \leq2e^{Ca_\infty}\int_{\mathbb R^N}F(X,T)\,dX,\ \ t\geq T,&\\
& \int_T^\infty\!\int_{\mathbb R^N}m(X,t,h(X,t))|\nabla_{X} h|^2\,dX\,dt
 \leq2e^{Ca_\infty}\int_{\mathbb R^N}F(X,T)\,dX.&
 \end{aligned}
\end{equation}
If \eqref{eq4.35} fails for some $H_0$, there exist
$d>0$, $t_j\to\infty$, $|{s}_j|\leq H_0$, and $y_j\in\mathbb R^{N-1}$
such that
\[
 h({s}_j,y_j,t_j)\geq d.
\]
Write $y_j=k_j+\widehat y_j$ and $t_j=n_jP+\vartheta_j$, 
where $k_j$ is a transverse period vector, $\widehat y_j\in\prod_{i=2}^N[0,L_i)$
and  $0\leq\vartheta_j<P$. After extraction,
\(({s}_j,\widehat y_j,\vartheta_j)\to({s}_*,y_*,\vartheta_*)\) as \(j\to+\infty\).
Define
\[
 h_j(X,\tau)=h(s,y+k_j,t_j+\tau), \ \ 
 v_j(X,\tau)=Z_{\sigma}((s+c_1(t_j+\tau),y+k_j),t_j+\tau).
\]
By periodicity and \eqref{eq4.28}, one has
\[\begin{gathered}
 (v_j)_\tau=\Delta_Xv_j+c_1(v_j)_{s}
 +f((s+c_1(\vartheta_j+\tau),y),v_j)+\sigma a'(t_j+\tau)J((s+c_1(\vartheta_j+\tau),y),v_j),\\
 v_j(X,\tau)=V^{\sigma(a(t_j+\tau)+h_j(X,\tau))}(X,\vartheta_j+\tau).
\end{gathered}\]
For every fixed $S>0$,
\[
 \sup_{|\tau|\leq S}|a(t_j+\tau)-a_\infty|
 \leq a_\infty-a(t_j-S)\to0,
 \ \ \sup_{|\tau|\leq S}|a'(t_j+\tau)|\to0.
\]
Since  $0\leq h_j\leq D$, the  Arzel\`a--Ascoli theorem and a diagonal
argument yield, along a subsequence, that 
\begin{equation*}
 h_j\to h_\infty
 \quad\text{locally uniformly on }\mathbb R^N\times\mathbb R,
 \ \  0\leq h_\infty\leq D.
\end{equation*}
Moreover, since $({s}_j,\widehat y_j)\to({s}_*,y_*)$,
\[
 h_\infty({s}_*,y_*,0)
 =\lim_{j\to\infty}h_j({s}_j,\widehat y_j,0)
 =\lim_{j\to\infty}h({s}_j,y_j,t_j)\geq d.
\]
Set $Q_H=(-H,H)\times B_H^{N-1}\times(-H,H)$.
By \eqref{eq4.34}, $m(X,t,h)\geq c_H>0$ on $|s|\leq H$,
uniformly in $y,t$. It follows from \eqref{eq4.37} that
\begin{equation*}
 \int_{Q_H}|\nabla_Xh_j|^2
 \leq c_H^{-1}\int_{t_j-H}^\infty\!\int_{\mathbb R^N}
   m(X,t,h)|\nabla_Xh|^2\,dX\,dt\longrightarrow0\ \ \text{as \(j\to+\infty\).}
\end{equation*}
For $\varphi\in C_c^\infty(Q_m)$,
\[
 \int_{Q_H} h_\infty\partial_{X_i}\varphi
 =-\lim_{j\to+\infty}\int_{Q_H}(\partial_{X_i}h_j)\varphi=0.
\]
Thus $\nabla_Xh_\infty=0$ in distributions. By  \eqref{eq4.31},  one has that 
\[
 \left|\int_{Q_H}h_j(-\varphi_\tau-\Delta_X\varphi)\right|
 \leq C\|\varphi\|_2\|\nabla_Xh_j\|_2
+C\|\varphi\|_\infty\|\nabla_Xh_j\|_2^2\longrightarrow0.
\]
Hence $\partial_\tau h_\infty-\Delta_Xh_\infty=0$. Together with
$\nabla_Xh_\infty=0$, this gives
\[
 \partial_\tau h_\infty=0,
 \ \  h_\infty\text{ is constant},\quad h_\infty\geq d>0.
\]
Let $H=1$. Then $F\geq c_Hh$ for $|s|\leq1$, and for every $L>0$,
\[
 c_H\int_{-1}^1\int_{B_L^{N-1}}h_j(s,y,0)\,dy\,ds
 \leq\int_{\mathbb R^N}F(X,t_j)\,dX\leq 2e^{Ca_\infty}\int_{\mathbb R^N}F(X,T)\,dX.
\]
Letting $j\to\infty$ yields
\[
 2c_H h_\infty|B_1^{N-1}|L^{N-1}\leq 2e^{Ca_\infty}\int_{\mathbb R^N}F(X,T)\,dX
 \quad\text{for every }L>0,
\]
which contradicts $N-1\geq1$ and $h_\infty>0$. This proves
\eqref{eq4.35}.

Given $\varepsilon\in(0,1)$, 
choose $H>c_1(D+a_\infty)$ such that, for all
$x\in\mathbb R^N$,
\[
 U_1(x,\xi)\leq\varepsilon\ \  \text{for \(\xi\geq H-c_1(D+a_\infty)\)},\ \
 1-U_1(x,\xi)\leq\varepsilon\ \ \text{for \(\xi\leq-H+c_1(D+a_\infty)\)}.
\]
Since $0\leq h\leq D$ and $0\leq a(t)\leq a_\infty$,
then $|\sigma(a(t)+h(X,t))|$ and $|\sigma a_\infty|$ are not more than $D+a_\infty$. 
Hence, for \(|r|\leq D+a_\infty\), one has that
\(V^r(X,t)=U_1((s+c_1t,y),s-c_1r)\)
belong to $[0,\varepsilon]$ when $s>H$,
and to $[1-\varepsilon,1]$ when $s<-H$. Consequently,
\[
 \sup_{\substack{|s|>H,\ y\in\mathbb R^{N-1}}}
 \left|V^{\sigma(a(t)+h(X,t))}(X,t)-V^{\sigma a_\infty}(X,t)\right|
 \leq\varepsilon.
\]
On $|s|\leq H$, the mean value theorem and
$0<\partial_rV^r=p^r\leq M$ give
\[
 \begin{aligned}
 \left|V^{\sigma(a(t)+h(X,t))}(X,t)-V^{\sigma a_\infty}(X,t)\right|
\leq M\bigl(h(X,t)+|a(t)-a_\infty|\bigr).
 \end{aligned}
\]
Combining the two regions yields
\[
 \begin{aligned}
 \sup_{X\in\mathbb R^N}
 \left|V^{\sigma(a(t)+h(X,t))}(X,t)-V^{\sigma a_\infty}(X,t)\right|
 \leq
 \max\left\{\varepsilon,\,
 M\left(
 \sup_{\substack{|s|\leq H\\y\in\mathbb R^{N-1}}}h(X,t)
 +|a(t)-a_\infty|
 \right)\right\}.
 \end{aligned}
\]
For this fixed $H$, \eqref{eq4.35} and
$a(t)\to a_\infty$ imply
\[
 \limsup_{t\to\infty}\sup_{X\in\mathbb R^N}
 \left|V^{\sigma(a(t)+h(X,t))}(X,t)-V^{\sigma a_\infty}(X,t)\right|
 \leq\varepsilon.
\]
Letting $\varepsilon\downarrow0$ and using
\eqref{eq4.30} proves \eqref{eq4.36}. 
This completes the proof.
\end{pr}
\vspace{0.3cm}

\begin{pr}[Proof of Theorem \ref{convergence}]
Fix $D$ as in Lemma \ref{lem-t17-trapping}. 
Recall that $\eta_0,\omega$ are fixed in \eqref{eq4.1},
$\rho$ is defined in \eqref{eq4.2}, and $\mu>0$
is chosen before \eqref{eq4.3}.
All these constants are independent of $\epsilon$.
Choose
\begin{equation}\label{eq4.38}
 0<\epsilon\leq {\eta_0}/2, \ \ 
  \rho\epsilon\left(\frac1\omega+\frac1{ \mu c_1}\right)\leq1.
\end{equation}
Fix $R>0$ such that $K\subset B_R$, where \( B_R:=\{x\in\mathbb R^N:|x|<R\}\)
and \(\partial B_R=\{x\in\mathbb R^N:|x|=R\}\).
Set $\mathcal D:=\mathbb R^N\setminus\overline B_R$.
By \eqref{eq4.6}, \eqref{eq4.7}, there is  $T>0$ so that, for every $t\geq T$,
\begin{equation}\label{eq4.39}
 \begin{gathered}
 W^{-D}(x,t)-\epsilon\leq u(x,t)\leq W^D(x,t)+\epsilon
 , \ \   x\in\overline\Omega,\\
 u(x,t),W(x,t)\geq1-\epsilon, \ \   x\in\partial B_R,
 \ \  e^{ \mu(R+ C_1+c_1(D+1))- \mu c_1T}\leq1.
 \end{gathered}
\end{equation}
By Theorem \ref{larget}, at $T$ there is a compact
$C_\epsilon\subset\overline\Omega$ such that
\begin{equation}\label{eq4.40}
 |u(x,T)-W(x,T)|\leq\epsilon
 , \ \   x\in\overline\Omega\setminus C_\epsilon.
\end{equation}
Since $C_\epsilon$ is compact, choose $L>0$ such that
\[
 C_\epsilon\subset
 \mathbb R\times\{y\in\mathbb R^{N-1}:|y|\leq L\}.
\]
Choose $\chi\in C_c^\infty(\mathbb R^{N-1})$ satisfying
\[
 0\leq\chi\leq1,\ \ 
 \chi(y)=1\quad\text{for }|y|\leq L,\ \ 
 \chi(y)=0\quad\text{for }|y|\geq2L.
\]
In particular, $\chi(y)=1$ whenever $(x_1,y)\in C_\epsilon$.
Set
\begin{equation*}
 \begin{gathered}
 h_0(s,y)=D\chi(y),\\
 a(t)= \rho\epsilon\left[
 \frac{1-e^{-\omega(t-T)}}{\omega}
 +\frac{e^{ \mu(R+ C_1+c_1(D+1))- \mu c_1T}}{ \mu c_1}
\bigl(1-e^{- \mu c_1(t-T)}\bigr)\right]
 , \ \   t\geq T.
 \end{gathered}
\end{equation*}
After some computations, it follows from \eqref{eq4.38} and \eqref{eq4.39} that 
\begin{equation}\label{eq4.41}
 \begin{gathered}
 a(T)=0,\quad
 a'(t)= \rho\epsilon
   \{e^{-\omega(t-T)}+e^{ \mu(R+ C_1+c_1(D+1))- \mu c_1t}\}\geq0,\ \   a'\to0 \ \text{as \(t\to+\infty\)},\\
 a''(t)=- \rho\epsilon
   \{\omega e^{-\omega(t-T)}
    + \mu c_1 e^{ \mu(R+ C_1+c_1(D+1))- \mu c_1t}\},\ \
\|a'\|_\infty+\|a''\|_\infty<\infty,\\
 0\leq a(t)\leq a_\infty
 = \rho\epsilon\left(\frac1\omega
  +\frac{e^{ \mu(R+ C_1+c_1(D+1))- \mu c_1T}}{ \mu c_1}\right) \leq  {\rho}\epsilon\left(\frac1\omega+\frac1{ \mu c_1}\right)\leq1.
 \end{gathered}
\end{equation}
Thus on can apply Lemmas \ref{lem-t17-auxiliary}--\ref{lem-t17-decay} to
$Z_-,Z_+$ with  $a$ and $h_0$ given above. In particular, it holds 
\begin{equation*}
 W^{-a-D}\leq Z_-\leq W^{-a}\leq W
 \leq W^a\leq Z_+\leq W^{a+D}, \ \ 
 Z_\pm(x,T)=W^{\pm D\chi(y)}(x,T).
\end{equation*}
For $x\in C_\epsilon$, $\chi(y)=1$, so \eqref{eq4.39}
gives $Z_-(x,T)-\epsilon\leq u(x,T)\leq Z_+(x,T)+\epsilon$.
For $x\notin C_\epsilon$, \eqref{eq4.40} and
$Z_-(x,T)=W^{-D\chi(y)}(x,T)
\leq W(x,T) \leq W^{D\chi(y)}(x,T)=Z_+(x,T)$ give the same inequalities. Hence
\begin{equation}\label{eq4.42}
 Z_-(x,T)-\epsilon\leq u(x,T)\leq Z_+(x,T)+\epsilon
 ,\ \ x\in\overline{\mathcal D}.
\end{equation}
On $\overline{\mathcal D}\times[T,\infty)$, define
\[
 \begin{aligned}
 w^-(x,t)&=Z_-(x,t)-\epsilon e^{-\omega(t-T)}
 -\epsilon e^{- \mu(|x|-R)},\\
 w^+(x,t)&=Z_+(x,t)+\epsilon e^{-\omega(t-T)}
 +\epsilon e^{- \mu(|x|-R)}.
 \end{aligned}
\]
Recall from \eqref{eq4.3} that
\( \zeta(x)=e^{-\mu(|x|-R)}\).
By \eqref{eq4.28}, one has that 
\[
 \begin{aligned}
 \Nop[w^-]=-a'J(x,Z_-)+\omega\epsilon e^{-\omega(t-T)}
  +\epsilon\Delta\zeta+f(x,Z_-)-f(x,w^-).
 \end{aligned}
\]
Since
\[
 |\Delta\zeta|\leq\omega\zeta, \ \ 
 0\leq\epsilon e^{-\omega(t-T)}+\epsilon\zeta
 \leq2\epsilon\leq {\eta_0},
\]
if $Z_-\leq {\eta_0}$ or $Z_-\geq1- \eta_0$, then
\eqref{eq4.5} gives
\[
 \Nop[w^-]\leq-a'J(x,Z_-)
 -3\omega\epsilon(e^{-\omega(t-T)}+\zeta)\leq0.
\]
If $ \eta_0\leq Z_-\leq1- \eta_0$, it then follows from   \(Z_-(x,t)=U_1\bigl(x,x_1-c_1(t-a(t)-h_-(X,t))\bigr)\)
 that there is  a constant \(C_1\) such that \(|x_1-c_1(t-a-h_-)|\leq  {C_1}\),
hence it can be inferred from  $0\leq h_-\leq D$ and \eqref{eq4.41} that 
\(|x|\geq x_1\geq c_1t-c_1(D+1)- C_1\), whence \(\zeta(x)
\leq e^{ \mu(R+ C_1+c_1(D+1))- \mu c_1t}.\)
Then by \eqref{eq4.25} one obtains that 
$J(x,Z_-)\geq  {\kappa}$, it then follows from   
\eqref{eq4.1}, \eqref{eq4.2} and \eqref{eq4.41} that 
\[
 \begin{aligned}
 \Nop[w^-]
 &\leq-a' \kappa+(L_f+\omega)\epsilon(e^{-\omega(t-T)}+\zeta)\\
 &\leq(- \rho {\kappa}+L_f+\omega)\epsilon
   \{e^{-\omega(t-T)}+e^{ \mu(R+ C_1+c_1(D+1))- \mu c_1t}\}=0.
 \end{aligned}
\]
On $\partial B_R$, it can be derived from \eqref{eq4.39} and $Z_-\leq1$ that
\[
 w^-(x,t)=Z_-(x,t)-\epsilon e^{-\omega(t-T)}-\epsilon
 \leq1-\epsilon\leq u(x,t).
\]
At $t=T$, it follows from  \eqref{eq4.42} that
\[
 w^-(x,T)=Z_-(x,T)-\epsilon-\epsilon e^{- \mu(|x|-R)}
 \leq Z_-(x,T)-\epsilon\leq u(x,T).
\]
A similar argument yields that 
\[
 \Nop[w^+]\geq0,\quad
 w^+\geq u\ \text{on }\partial B_R\times[T,\infty),\quad
 w^+(x,T)\geq u(x,T).
\]
It then follows from the comparison principle  that
\begin{equation}\label{eq4.43}
 \begin{aligned}
 Z_--\epsilon e^{-\omega(t-T)}-\epsilon e^{- \mu(|x|-R)}
 \leq u\leq Z_++\epsilon e^{-\omega(t-T)}+\epsilon e^{- \mu(|x|-R)}\ \ \text{in  $\overline{\mathcal D}\times[T,\infty)$.}
 \end{aligned}
\end{equation}
For fixed $\epsilon$, Lemma \ref{lem-t17-decay} gives
\[
 \max_{\sigma=\pm1}
 \|Z_{\sigma}(\cdot,t)-W^{\sigma a_\infty}(\cdot,t)\|_{L^\infty(\mathbb R^N)}
 \longrightarrow0\ \ \text{as \( t\to+\infty\),}
\]
where \(W^{\sigma a_\infty}(x,t)=W(x,t+\sigma a_\infty)\).
The mean value theorem and \eqref{eq4.2} give
\[
 \|W^{\pm a_\infty}(\cdot,t)-W(\cdot,t)\|_{L^\infty(\mathbb R^N)}
 \leq Ma_\infty.
\]
By \eqref{eq4.43}, in $\overline{\mathcal D}\times[T,\infty)$, there hold
\[
 \begin{aligned}
 u-W&\leq |Z_+-W^{a_\infty}|+Ma_\infty
+\epsilon e^{-\omega(t-T)}+\epsilon,\\
 W-u&\leq |Z_--W^{-a_\infty}|+Ma_\infty
+\epsilon e^{-\omega(t-T)}+\epsilon.
 \end{aligned}
\]
Hence 
\[
 \limsup_{t\to\infty}
 \|u(\cdot,t)-W(\cdot,t)\|_{L^\infty(\mathcal D)}
 \leq Ma_\infty+\epsilon
 \leq\left[M \rho\left(\frac1\omega+\frac1{ \mu c_1}\right)+1\right]\epsilon.
\]
In addition, by \eqref{eq4.7}, one has 
\[
 \sup_{x\in\overline\Omega\cap\overline B_R}|u(x,t)-1|
 \longrightarrow0
 \quad\text{as }t\to\infty.
\]
Since $x_1-c_1t\leq R-c_1t$ on $\overline B_R$, then 
\[
 \begin{aligned}
 \sup_{x\in\overline\Omega\cap\overline B_R}|W(x,t)-1|
 &\leq
 \sup_{\substack{x\in\mathbb R^N\\ \xi\leq R-c_1t}}
 |U_1(x,\xi)-1|
 \longrightarrow0
 \ \ \text{as }t\to\infty.
 \end{aligned}
\]
Consequently,
\[
 \begin{aligned}
 \sup_{x\in\overline\Omega\cap\overline B_R}|u(x,t)-W(x,t)|
 \leq
 \sup_{x\in\overline\Omega\cap\overline B_R}|u(x,t)-1|
+ \sup_{x\in\overline\Omega\cap\overline B_R}|W(x,t)-1|
 \longrightarrow0
 \ \ \text{as }t\to\infty.
 \end{aligned}
\]
Thus
\[
 0\leq\limsup_{t\to\infty}
 \|u(\cdot,t)-W(\cdot,t)\|_{L^\infty(\overline\Omega)}
 \leq\left[M \rho\left(\frac1\omega+\frac1{ \mu c_1}\right)+1\right]\epsilon.
\]
The constants on the right are independent of $\epsilon$.
Letting $\epsilon\downarrow0$ one arrives at
\[\lim_{t\to\infty}\sup_{x\in\overline\Omega}
 |u(x,t)-U_1(x,x_1-c_1t)|=0.
\]
This completes the proof.
\end{pr}

\section*{Acknowledgments}
This work was supported by NSF of China (12671234, 12171120) and by the Fundamental Research Funds for the Central 
Universities.

%
%\section*{Date availability statements}
%We do not analyse or generate any datasets, because our work proceeds within a theoretical and mathematical approach.
%
%\section*{Conflict of interest}
%There is no conflict of interest to declare.
%
%\section*{Declaration of generative AI use}
%The generative AI was not used in this article.

\end{document}